\documentclass[12pt]{article}
\usepackage{amsmath,amscd,amsthm,amsfonts,amssymb}
\usepackage{array}
\usepackage{graphics}
\usepackage{amscd}
\usepackage{graphicx,epsfig}
\usepackage{xcolor}
\usepackage{fancyhdr}
\def\goth#1{\mathfrak{#1}}

\newtheorem{thm}{Theorem}

\newtheorem{ex}{Example}
\newtheorem{cor}{Corollary}
\newtheorem{lem}{Lemma}

\newtheorem{prop}{Proposition}

 \definecolor{marronf}{RGB}{97, 22, 22}
  \definecolor{bleuf}{RGB}{8, 8, 172}
   \definecolor{vertf}{RGB}{25, 122, 10}
    \definecolor{orangef}{RGB}{201, 77, 16}
     \definecolor{rougef}{RGB}{182, 13, 7}
\definecolor{olive}{rgb}{0.3, 0.4, .1}
\definecolor{fore}{RGB}{249,242,215}
\definecolor{back}{RGB}{51,51,51}
\definecolor{title}{RGB}{255,0,90}
\definecolor{dgreen}{rgb}{0.,0.6,0.}
\definecolor{gold}{rgb}{1.,0.84,0.}
\definecolor{JungleGreen}{cmyk}{0.99,0,0.52,0}
\definecolor{BlueGreen}{cmyk}{0.85,0,0.33,0}
\definecolor{RawSienna}{cmyk}{0,0.72,1,0.45}
\definecolor{Magenta}{cmyk}{0,1,0,0}
\newcommand{\ndiv}{\hspace{-2pt}\not\hspace{-2pt}|\hspace{2pt}}
\begin{document}
\title{ Lissajous and Fourier  knots  }
\author{Marc Soret and Marina Ville }
\date{ }
\maketitle
\begin{abstract}
We prove that any knot of $\mathbb{R}^3$ is isotopic to a Fourier knot of type $(1,1,2)$ obtained by 
deformation of a Lissajous knot.  \end{abstract}

\section{Introduction}
Fourier knots are   closed embedded curves  whose coordinate functions are finite Fourier sums. 
   Lissajous knots  (defined in    \cite{BHJS})  are 
   the  simplest examples : 
   each coordinate function  consists in  only one term. Lissajous knots  are Fourier knots of type (1,1,1) (cf.   for example   \cite{K}, \cite{La2},\cite{T}). \\
 \indent   Surprisingly - at least at first sight - not  every  isotopy class of knots 
  can be represented by Lissajous knots. 
  Indeed Lissajous knots   are isotopic to their mirror image;  in particular  a nontrivial 
  torus knot cannot be  isotopic to  a Lissajous knot.\par
 Let us first  recall how to construct knots from a knot shadow. \par
A knot shadow is  a generic projection  of a knot on a  plane. 
It is a  closed  planar curve 
 with nodes, i.e.  double points.
Conversely   given a  shadow -  i.e. an oriented  planar closed curve with double points :
 $$D:  t \mapsto \gamma(t) =\left(x\left(t\right),y\left(t\right)\right)\in \mathbb{R}^2,$$
 we can  construct a  knot  in   $\mathbb{R}^3$  
 by   defining   a  height  function  ${\bf z}$  which has the right values at each node of the shadow.
The knot is then defined by: 
 $${\bf  K}:  t \mapsto \left(x\left(t\right),y\left(t\right), {\bf  z\left(t\right)}\right)\in \mathbb{R}^3.$$
 % \begin{figure}[!t]
 %\fbox{  \includegraphics[width= 5 cm, height = 4 cm]{image_construire.png}}
 % \end{figure} 
Thus  a Lissajous knot   projects onto  a Lissajous shadow  ($x(t) $ and $y(t)$ are cosine functions). 
Although we cannot represent any knot by a Lissajous knot, still any knot $K$ 
is isotopic to a knot which projects onto  a Lissajous shadow (\cite{La}).
 In other words one can choose $D$ to be a Lissajous shadow   
  but one cannot  always choose  {\bf $ z(t) = \cos( pt + \phi)$} as the height function.
However, this is possible  if the height function $z$  is  a Fourier sum of a non prescribed finite number
of terms (cf.  for instance \cite{K}).\\
  It was conjectured  in \cite{La2}  that any knot can be presented by a Lissajous diagram
  with a height function consisting of  a Fourier sum with a fixed number of terms 
 or  even,  as  it was  suggested experimentally in  \cite{BDHZ},  with a height function consisting of  only two terms:
namely 
  \begin{thm}\label{lissajous} Any knot in $\mathbb{R}^3$   is isotopic to a Fourier knot of type (1,1,2).
   \end{thm}
 %The main idea is  perturb slightly  a lissajous curve by  adding a  small term of type  $\epsilon\cos(n t+\phi )$ to ,say the $y$ coordinate 
%  of the Lissajous  parametrization. The nodes parameters of the  Theorem 2 is more a  ``middle range" result about minimal knots in the terminology of Rudolph, but still sheds new light on the 
%local nature of minimal knots. \\
  
The technique of the  proof of theorem \ref{lissajous}      is inspired by  a paper of \cite{KP} and uses 
Kronecker's theorem. Another key property of number theory that will also   be used  is  related to  the fact that the only rational  values 
of $\sin \frac{2\pi p}{q} $- for integers $p$ and $q$- are only $0, \pm 1, \pm 1/2$.\\
 The main  idea  goes as follows:
given a knot $K$ in a given isotopy class of knots, we show that it can be presented by   a Lissajous diagram   with prime frequencies.
We then deform the shadow so that nodes are in ``general position",  i.e.   their parameters are 
rationally linearly independent. This will be derived  from
the fact that the nodal curve of the deformation  is skewed in the parameter space of all nodes. We  then  find an integer   $n$ and   real number $\phi$ -small-  such that 
 the  height function   $ h(t) := \cos\left[2n\pi\left(t+\phi\right)\right]$ 
 has  the right intersection signs for every  two  values of the parameter $t$ corresponding to  a node  of the Lissajous shadow. 
 The height function   then defines  a third coordinate function 
 which together with the deformed Lissajous shadow  coordinates defines a curve of $\mathbb{R}^3$ in the given isotopy class. \par
 The paper is organized as follows :
 In Section  2 and 3  we fix notations and terminology. We compute  the  node parameters 
  of deformations of Lissajous planar curves in section 4; we also show that any knot admits a Lissajous diagram
  with prime frequencies.\\
  In Section 5, we show that the  nodal curves of our deformations are skewed :
 for an appropriate choice of deformation parameters and  a simple choice of phases, we can compute the determinant 
 of the $k$-th derivatives $k=1,\cdots, n_D$    of the nodal curve
   position vector in $\mathbb{R}^{n_D}$  and show that  it is not zero. 
    In Section \ref{sub3}, the  determinant is expressed as a product of factors,  
 each of which is proved to be different from zero.
We also show that the skewness of the nodal curve implies  that    nodes parameters  of the diagram  are rationally linear independent
for a dense subset of the deformation parameter.
The last section is devoted to the proof of the theorem  :  Kronecker's theorem applied to  $\mathbb{Q}$-linear  independence  
of nodal coordinates allows us  to  choose, as  height function 
of  any knot represented by an appropriate  Lissajous  shadow,  a cosine function with appropriate frequency. \\
{\it Acknowledgement} : the starting point of this paper was a fruitful meeting with P-V Koseleff whom we thank 
heartfully.
 
 \tableofcontents
\newpage
\section{Terminology}
\subsection{Knots in $\mathbb{R}^3$}\label{term1}
 We will consider knots in $\mathbb{R}^3$  as  presented  by   a smooth embedding of the circle:
 \begin{equation}\label{shadow}
 \left( 
 \begin{array}{cc}
 \gamma :\mathbb{S}^1 = \mathbb{R} / \mathbb{Z} &\longrightarrow \mathbb{R}^3 =\mathbb{R}^2\times \mathbb{R} \\
 t   &\mapsto  \left( \gamma_D(t), z\left(t\right)\right)
 \end{array}
 \right)
 \end{equation} 
where  the {\it height function } $z$ is a  smooth function $\mathbb{S}^1 \longrightarrow \mathbb{R}$.  The
planar  curve $\gamma_D$  is a {\it  knot shadow}, 
if it is  furthermore a smooth  immersion  such that  the  self-intersections of $\gamma_D$  are   transverse double points also called {\it nodes}.
Hence for any node $P\in  \gamma_D(\mathbb{S}^1)$ there  exists a pair of distinct real numbers $s,t \in [0,1[$ such that 
$\gamma_D(s)=\gamma_D(t)=P$.  We choose for each node one of these two parameters and   denote it  by $t_P\in [0,1[$. 
It will be convenient to  denote   the set of nodes by an ordered  set  $J$ of   points and  define accordingly the  {\it   nodal vector} $\eta_0$  of $\gamma_D$    by  the $n_D$-uplet   $( t_P)_{P\in J}$  $[0,1[^{n_D}$,  where  $n_D$ is the number of nodes  of 
curve $\gamma_D$. \\
Conversely  given a knot shadow $\gamma_D$, and a height function $z(t)$, we can define a knot  as in expression \eqref{shadow}
if, for each node $P$ of the shadow and for each corresponding pair of  parameters $\{s,t\}$ of $P$ , $z(s)-z(t) \not = 0$.  In fact 
the knot thus obtained is entirely defined by curve $\gamma_D$ and the  data  $ sign(z(s)-z(t) )$ for each node.
In other words a knot is defined by a shadow $D$  and a height function $z$, where each  node of $D$ is parametrized by a pair of  real numbers $\{s,t\}$ such that $\gamma_D(s)=\gamma_D(t)$ and the 
  sign  of  $ z(t) -z(s) $   defines  which strand of $\gamma$ lies  above which one at a crossing point-or node-  of $\gamma_D$.\\
% \begin{figure}[!h]
 %\fbox{  \includegraphics[width= 4 cm, height = 4 cm]{image_shadowme}}
% $\includegraphics[width= 2.5 cm, height = 1 cm]{image_crossing.png} }
  %\end{figure} 
 
 \subsection{Deformation of a knot  shadow and nodal curve}\label{nodalcurve}
A  {\it  deformation} of a shadow $\gamma_D$  is a  smooth   family of curve shadows 
$\{\gamma_{D,\epsilon}\}_{\epsilon \in I}$ 
where $\epsilon $ lies in the {\it interval of deformation} $ I= [0,\epsilon_0] $ 
and such that $\gamma_{D,0}=\gamma_D$.  We associate  to a shadow deformation   the {\it nodal curve }
$\eta : [0, \epsilon_0]  \longrightarrow [0,1]^{n_D} $ 
  defined as follows.  For each  node $P_i$, 
$ i = 1, \cdots, n_D$, of the  deformed shadow and  ordered in a certain way, 
there is a pair of two parameters $\{s_i(\epsilon),t_i(\epsilon)\}$ such that 
$\gamma_D(s_i) = \gamma_D(t_i) = P_i$. We choose one of the two parameters, say $t_i$ and the nodal curve is then 
defined by $\eta(\epsilon)   =\left(  t_j \left(\epsilon\right)\right)_{j= 1,\cdots, n_D }.  $ 
 \section{Fourier  Knots}
 We give a quick reminder of   some geometric knot presentations that are of some  interest in our topic:  Fourier knots, 
 Lissajous knots and  torus knots  that all belong to next family.\\
  A Fourier knot -of type $(m,n,p)$- is a  parametrized  curve of  $\mathbb{R}^3$ defined  by :
  \begin{equation}\label{fourier}
 F_{m,n,p}: \left(
  \begin{array}{cc}
[0,1[  &\longrightarrow\mathbb{R}^3 \\
    t    & \mapsto  F_{m,n,p}(t)      
   \end{array}
    \right)
 \end{equation}
 with 
  $$F_{m,n,p}(t)      = \left(  \sum_{i=1}^m M_i \cos(2\pi m_i t +\phi_i),
 \sum_{i=1}^n N_i \cos(2\pi n_i t +\psi_i),  \sum_{i=1}^p P_i \cos(2\pi p_i t +\tau_i)\right)$$
 for some  $\phi_i,M_i,N_i ,\psi_i,, \tau_i, P_i\in \mathbb{R}$ and $ m_i,n_i,p_i \in \mathbb{N}$. \\
 It was proved in \cite{La2}  that    any   knot $K$  is isotopic to  a
 curve of type $(1,1, n_K)$  where $n_K$ depends on the knot.
 \subsection{Lissajous  knots}
   Fourier knots of type  $(1,1,1)$  are  also called Lissajous knots and have been extensively studied in \cite{BHJS}). 
\begin{equation}
 L(n_1,n_2,n_3,\phi_1,\phi_2) :\left( 
\begin{array}{ll}
  [0,1] &\longrightarrow \mathbb{R}^3\\
  t &\mapsto   \left(     \cos 2\pi n_1 t ,
 \cos2\pi n_2 \left(t +\phi_1\right),\cos2\pi n_3\left( t +\phi_2\right)   \right)
  \end{array}
  \right)
 \end{equation}

 \begin{figure}[!h]
\fbox{ \includegraphics[width= 4 cm, height = 4 cm]{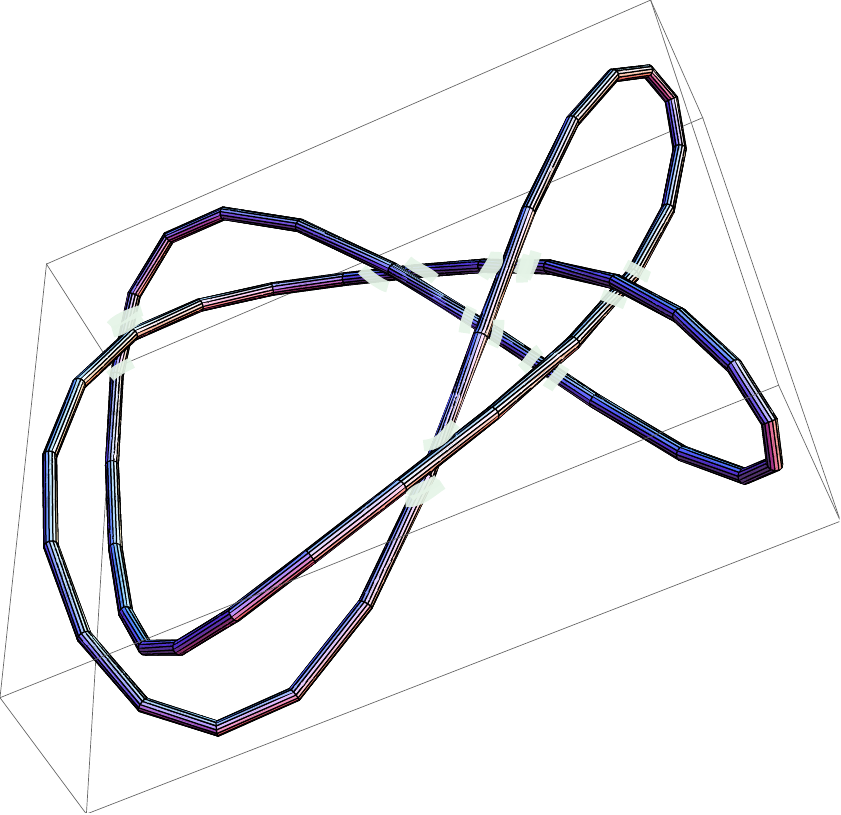}
\includegraphics[width= 4 cm, height = 4 cm]{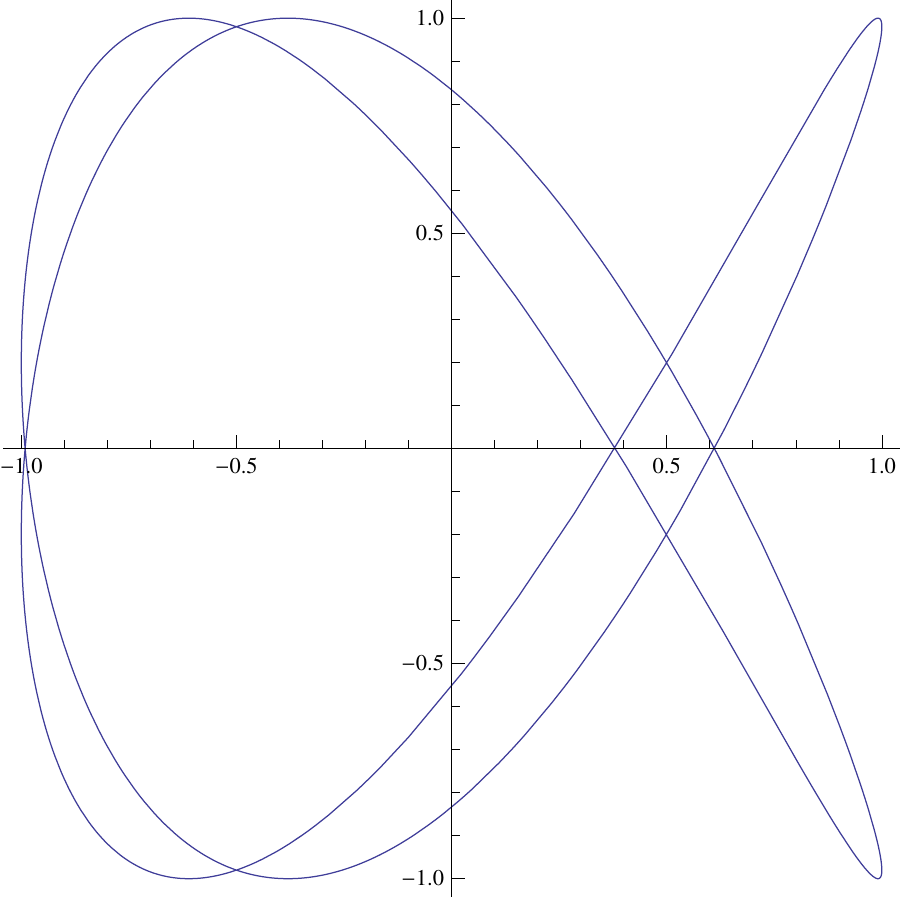},
\includegraphics[width= 4 cm, height = 4 cm]{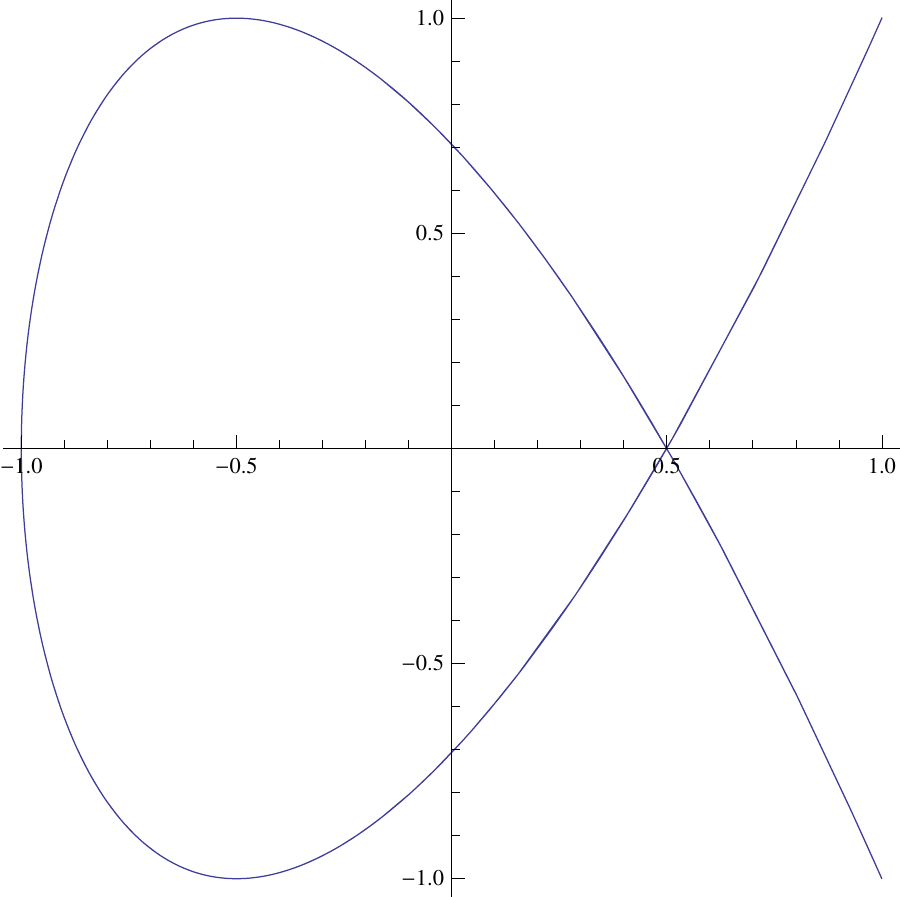}
 }
\caption{ knot $L(2,3,5,0,.2,2)=6_1$, shadow $L(2,3,.2)$, and associated Chebyshev figure $C(2,3)$}
\end{figure}
It is noteworthy to recall that Lissajous knots   are topologically equivalent to  closed billiard trajectories in a cube. 
(cf. \cite{JP}) 
% \begin{figure}[!t]
%\fbox{ \includegraphics[width= 3 cm, height = 2 cm]{liss4-5-gen}, 
%\includegraphics[width= 3 cm, height = 2 cm]{liss4-5-cheb}}
%\caption{L(}
%\end{figure} 
These Lissajous knots project horizontally   on planar Lissajous  curves  of type  $L(n_1,n_2, \phi_1)$. \\
When the  phase $\phi_1$  is zero,  the Lissajous figure 
degenerates into a 2-1 curve $C(n_1,n_2)$ which is a subset of an algebraic  Chebyshev  open curve  as defined in \cite{KP}: 
$$T(n_1,n_2) = \{ (x,y) :  T_{n_1}(x) = T_{n_2}(y) = 0\},$$  where 
$T_n$  is the Chebyshev polynomial of degree $n$ (see fig. 1). \\
 Although not all knot are Lissajous knots, \cite{KP}  showed that, for suitable numbers $n_1,n_2,n_3$ and phase $\phi$,
 any knot is isotopic  a {\it Chebyshev knot} obtained by Alexandrov compactification of  a curve  
 $T(n_1,n_2,n_3,\phi) : \mathbb{S}^1 \longrightarrow \mathbb{R}^3$  with  
  $  T(n_1,n_2,n_3,\phi)(t) =  ( T_{n_1}(t), T_{n_2}(t),T_{n_3}(t+\phi )).$
 \subsection{Torus  knots}\label{torus}
Another key family of Fourier knots 
are torus knots -neither of which can be  isotopic to a Lissajous knot (except the trivial ones)!
 A Torus knot $T(p,q)$ is originally  defined as  an embedding  of $\mathbb{S}^1$  onto  a torus as a Fourier knot  of type $(3,3,1)$  : 
 $$ T(p,q)(t)  = \left(  \cos 2\pi qt\left( 1 +\frac{1}{2} \cos2\pi p t\right) ,
 \sin  2\pi qt \left( 1 +\frac{1}{2} \cos 2\pi p t  \right),
   \sin  2\pi  pt \right).$$
  %   \begin{figure}[!h]
%\fbox{ \includegraphics[width= 6 cm, height = 4 cm]{image_torus}, 
%\includegraphics[width= 3 cm, height = 2 cm]{image_tore2}}
%\caption{T( 2, 3)}
%\end{figure} 
But it was shown   in \cite{H} that   torus knots  $T(p,q)$'s are  isotopic to  
 Fourier knots of type  $(1,1,2)$ :
 $$ T(p,q)(t)  = \left(  \cos 2\pi pt  ,
  \cos 2\pi q \left(t +  \frac{1}{4p}\right)  ,  \cos 2\pi \left(pt +\frac{1}{4}\right)    + \cos 2\pi \left(\left(q-p\right)t +\frac{1}{4p}\right) \right)$$

 \section{Nodes and deformations of  knot shadows}
We  first recall some basic properties of   Lissajous planar curves  and  compute  the  positions of its nodes.\\
In the second part of the section we add   a  small perturbation term to one of the two coordinates;
  we  define  a family  of perturbed  Lissajous planar curves  which are very close to the
 original Lissajous figure. We then describe the positions of the nodes as functions of the  deformation parameter $\epsilon$.
 
 \subsection{Lissajous Figures, nodes and symmetries}
  
We need to give some precisions on  the Lissajous curves :
\begin{equation}
L(n_1,n_2,\phi ) : \left( 
\begin{array}{ll}
[0,1] &\longrightarrow \mathbb{R}^2 \\
t  &\mapsto   \left(     \cos\left(2\pi n_1 t  \right),
 \cos\left(2\pi n_2 t +\phi\right)  \right)
 \end{array}\right)\end{equation}
 
 We  will    suppose that  $n_1 $ and $n_2$  are coprime  (otherwise $L$ is not 1-1). 
 We can also suppose  for our purpose that the phase  $\phi $ is {\it a small positive irrational number}.
 Our first goal is to find the expression of the ordered pairs  of parameters  $(t_i,s_i), i = 1 ,\cdots N$ corresponding to the double points $P_i, i=1,\cdots, N$
 of a Lissajous figure.
 \begin{figure}[!h]
\fbox{ \includegraphics[width= 4 cm, height = 4 cm]{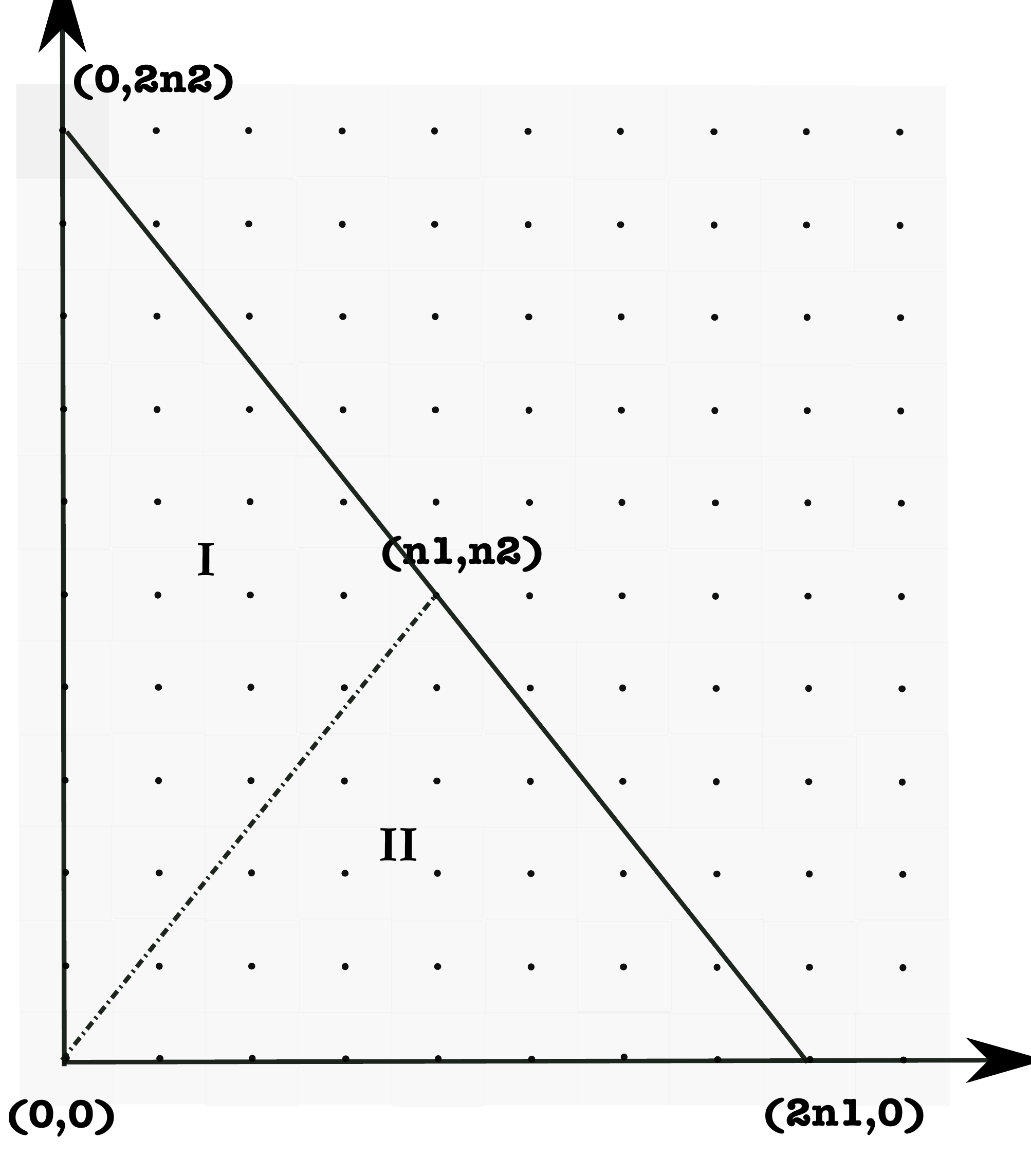}
\includegraphics[width= 4 cm, height = 4 cm]{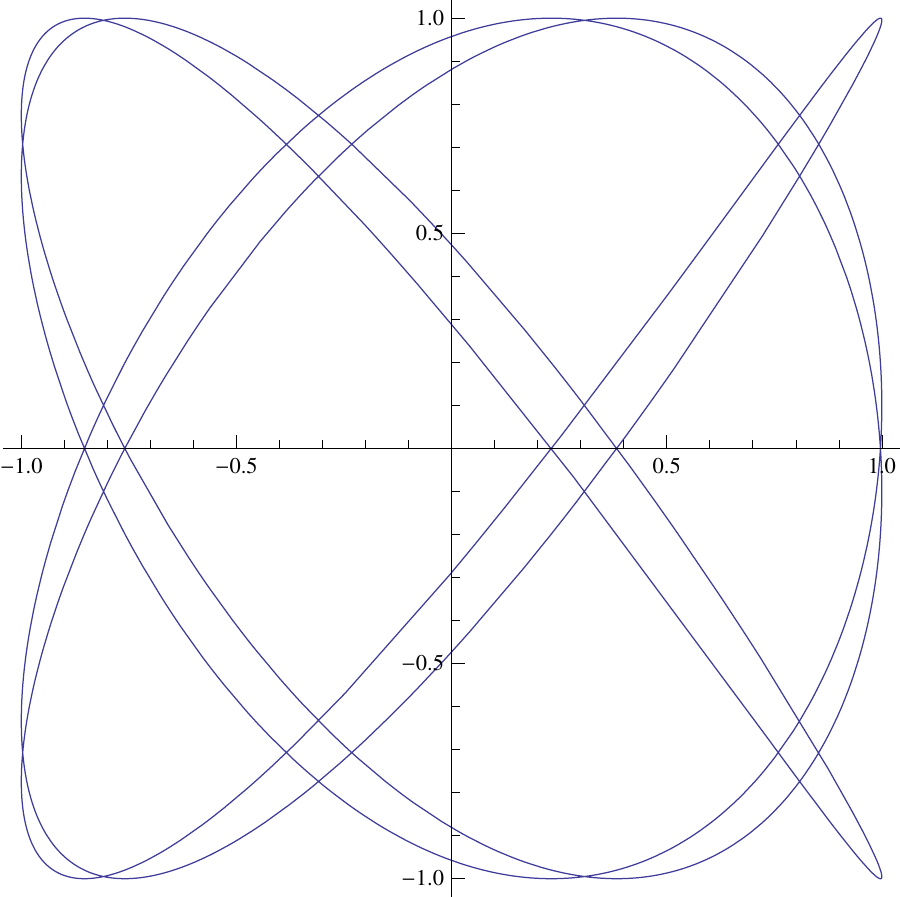},
\includegraphics[width= 4 cm, height = 4 cm]{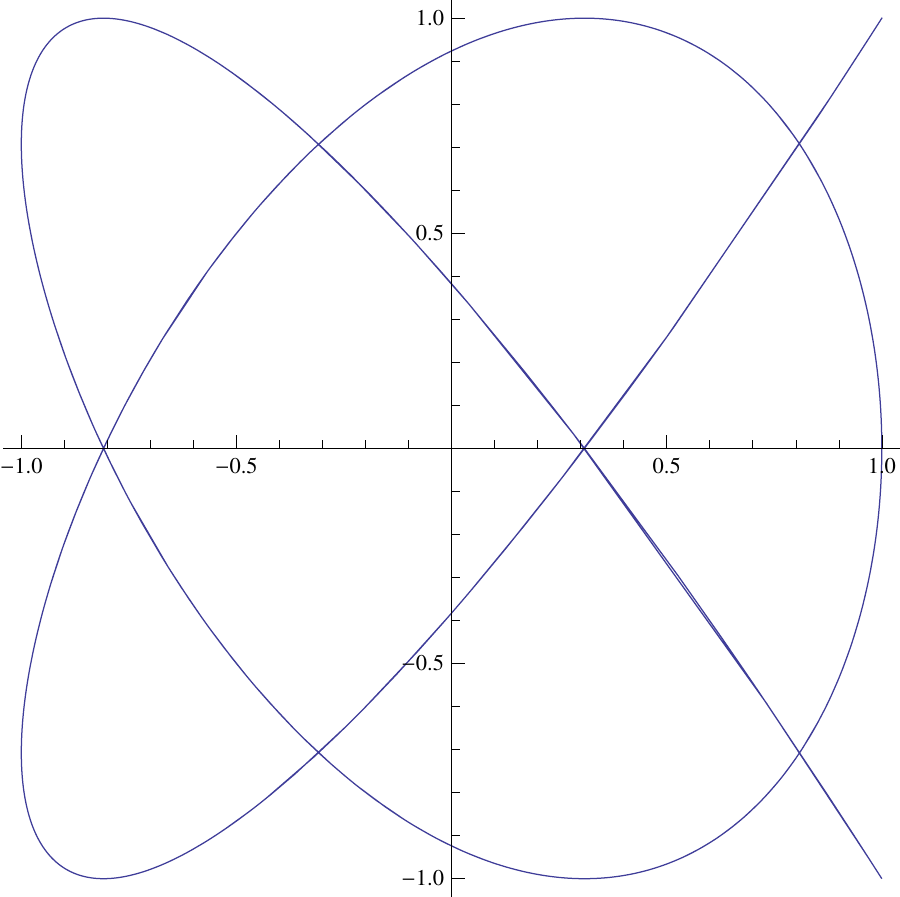}
 }
\caption{ 31 Nodes parametrized by integer points in  triangle $\Delta$  of Lissajous  $L(4,5,.1)$,   and  associated Chebyshev  $C(4,5)$ }
\end{figure}
\begin{figure}[!h]
\fbox{ \includegraphics[width= 4 cm, height = 4 cm]{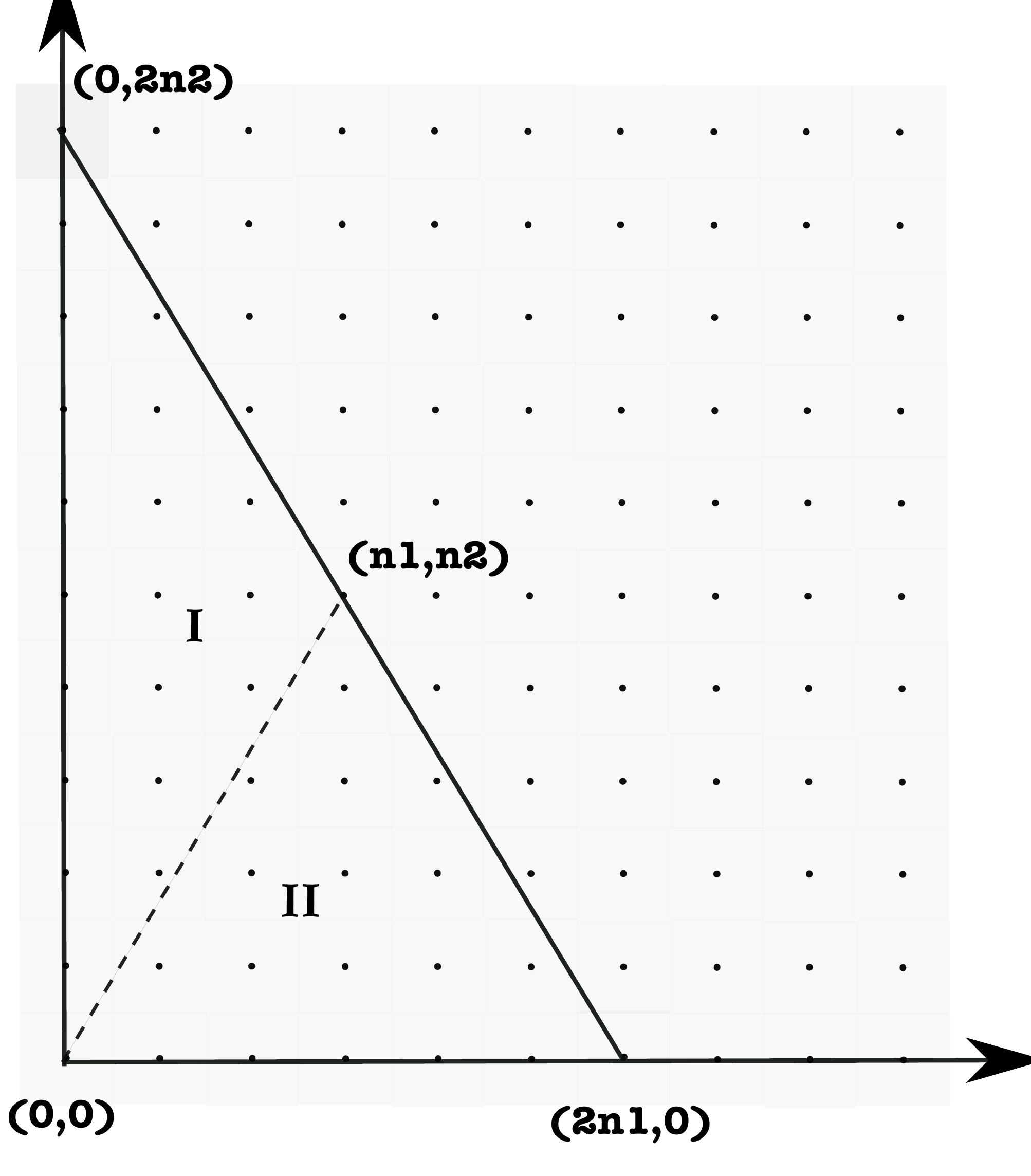}
\includegraphics[width= 4 cm, height = 4 cm]{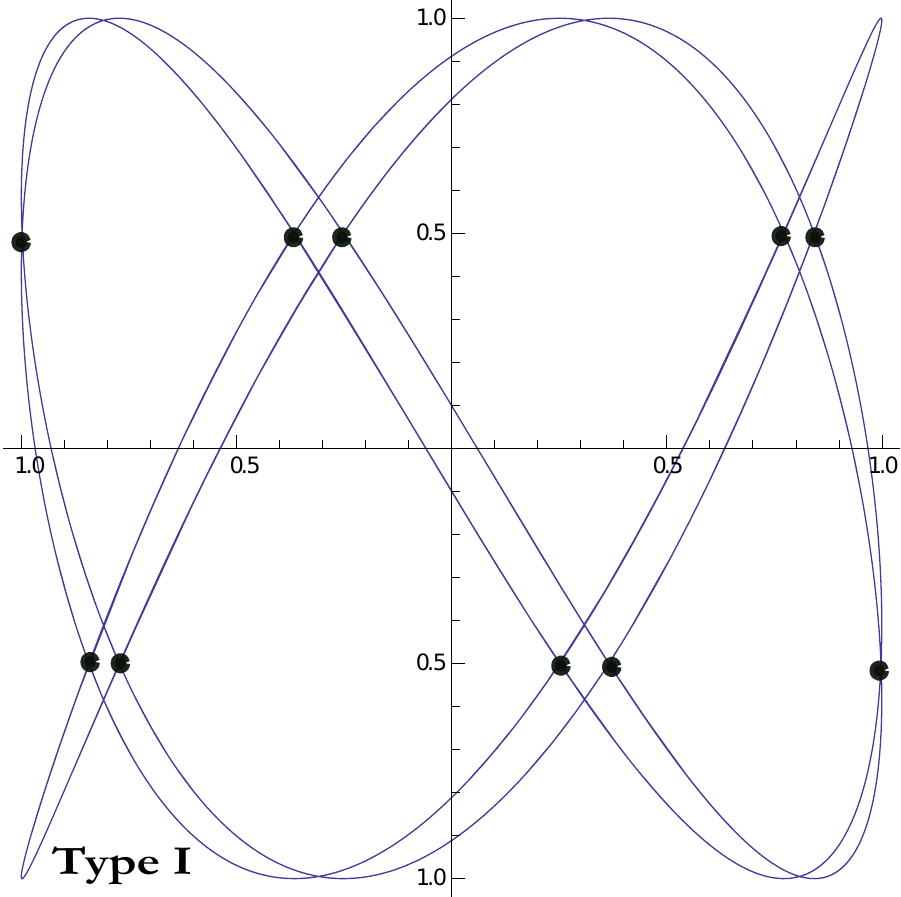},
\includegraphics[width= 4 cm, height = 4 cm]{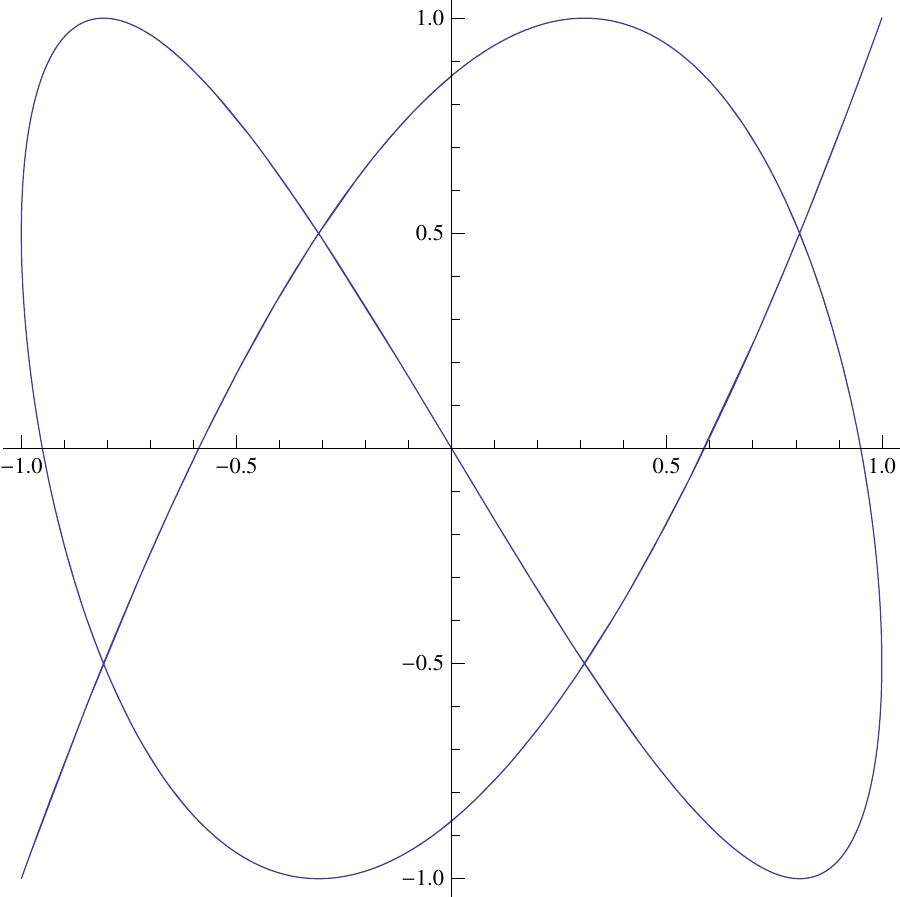}
 }
\caption{ 22 nodes parametrized by integer points    of  Lissajous  curve  $L(3,5,.1)$ with 10 nodes of type I,  and  associated Chebyshev  $C(3,5)$}
\end{figure}
In the  degenerate case where the phase $\phi$ is zero we obtain a $2-1$ curve (see figure 2 and 3). 
This curve is a subset of a Chebyshev curve defined in [KP] and  which is  an  open algebraic curve.
\subsection{Nodes parameters  of Lissajous planar curves}
We need to compute  the node parameters of a planar Lissajous curve. 
Such  a parametrization  was already published (cf. for instance  \cite{BDHZ} or \cite{JP}). But we derive it again
to get 
a geometric representation of these nodes as a set of integer points that lie in a  straight rectangle:  we will  use  symmetries 
of this set for our computations. \\
We may notice first that the number of nodes of a Lissajous figure is easily 
deduced from the number of nodes $ n_{C(n_1,n_2)}:= \frac{(n_1 -1)(n_2 -1)}{2}$ of its
  associated Chebyshev curve  (cf. for example \cite{Pe}).  As  the phase $\phi$ becomes positive,  each   node of 
 the Chebyshev curve blows up into  four nodes. Moreover  each pair of  maxima or minima of 
the coordinates functions  $x(t), y(t)$   give  also rise to a node (except for $t=0,\pm 1,\pm 1/2$) ; 
this produces  $n_{L(n_1,n_2,\phi)} = 4\frac{(n_1 -1)(n_2 -1)}{2} + n_1 +n_2 -2 = 2n_1n_2-n_1-n_2$.\\
  More precisely,  given a  node  $ P$, let us find  the double parameters $(t_P,s_P)$ such that 
$L(t_P) = L(s_P)=P$.\\
Equality for  the first coordinate  yields : 
$$\cos   2\pi n_1 t    = \cos    2\pi n_1 s   $$
iff
\begin{equation}\label{eq1}
s = \sigma t + \frac{k}{n_1},\quad  \sigma = \pm 1 \quad k\in\mathbb{Z}
\end{equation}
Equality for  the second coordinate  yields  : 
\begin{equation}\label{eq22}
\cos  \left[ 2\pi n_2  \left( t +\phi \right) \right]  = \cos  \left[  2\pi n_2 \left( s +\phi \right) \right] 
\end{equation}
Plug  equality \eqref{eq1} into \eqref{eq22} :
 \begin{equation}\label{eq2}
\cos     2\pi n_2  \left( t +\phi \right)   = \cos     2\pi n_2 \left( \sigma t + \frac{k}{n_1} +\phi \right)  
\end{equation}
and apply again \eqref{eq1} :
 \begin{equation}\label{eq3}
\sigma t + \frac{k}{n_1} +\phi  = \sigma' \left(     t +\phi             \right) + \frac{l}{n_2} \quad  \sigma' =\pm 1
 \quad k\in\mathbb{Z}\end{equation}
 We get two sets of solutions  according to the choice of $\sigma$ and $\sigma'$ :
 \begin{enumerate}
 \item{\bf Type I.}
  If  $\sigma =1$  and   $\sigma' =-1$  then : 
\begin{equation}\label{eq3}
t =-\phi +  \frac{1}{2}\left( \frac{l}{n_2}-\frac{k}{n_1} \right) 
\end{equation}
Hence 
\begin{equation}\label{eq3}
s = t +\frac{k}{n_1} =-\phi +  \frac{1}{2}\left( \frac{l}{n_2}+\frac{k}{n_1} \right) 
\end{equation}
 Since $s,t \in [0,1[$,  the  integer points  $(k,l)$ necessarily  lie  in  a parallelogram $\mathcal{P}$ which is defined by conditions :
\begin{equation}
\left\{  
\begin{array}{ccc}
2\phi n_1n_2\leq & n_1 l - n_2 k < & (1+\phi ) 2n_1n_2 \\
2\phi n_1n_2\leq &  n_1 l +  n_2 k < &(1+\phi ) 2n_1n_2 
\end{array}
\right.
\end{equation}
As phase  $\phi$ is small, parallelogram $\mathcal{P}$  
  is obtained by  a small  translation of   the parallelogram with  vertices $(0,0), (n_1,n_2),  (-n_1,n_2), (0,2n_2)$.
Since $n_1$ and $n_2$ are coprime   there are only the 4 integer points vertices on the boundary of
the translated parallelogram.\\

 \item{\bf Type II.}
In the second case,    $\sigma =-1$, $\sigma' =1$  and  we obtain : 
\begin{equation}\label{eq3}
t = \frac{1}{2}\left(\frac{k}{n_1} -\frac{l}{n_2} \right) 
\end{equation}
Hence 
\begin{equation}\label{eq3}
s =   \frac{1}{2}\left( \frac{l}{n_2}+\frac{k}{n_1} \right) 
\end{equation}
And  conditions $s,t\in [0,1[$ yields - we recall that $\phi$ is very small-
\begin{equation}
\left\{  
\begin{array}{ccc}
0\leq &  n_2 k- n_1 l  < &  2n_1n_2 \\
0\leq & n_2 k + n_1 l< &    2n_1n_2 
\end{array}
\right.
\end{equation}
 Solutions are given by integer points  $(k,l)\in \mathcal{P}'$ which is a small translation of 
the  parallelogram   with  vertices $(0,0), (n_1,n_2),  (n_1,-n_2), (2n_1,0)$.
\end{enumerate}
%In both cases, since  $n_1$ and $n_2$ are relatively prime hence for each $j$ such that  $ 0\leq j < 2n_1n_2$,   we can find 
%$(k,l)$  such that $j$ is a linear combination of $n_1$ and $n_2$ which gives  $2n_1 n_2$ possible integer points.\\
For each solution found, one should check that  $s\not = t$ 
and  also that   pairs $\{ t,s\}$ corresponding to the same node  are not   counted twice. \\
For each pair of solutions $\{ s,t\}$, condition $s\not = t$  
implies  $k \not = 0$ for the  first type
 and  $l \not = 0$ in the second case; furthermore for points of type I  the transformation : $k\mapsto -k$ permutes $t$ and $s$ 
  which are parameters that  correspond to 
 the same node. Hence we can reduce the  parameter set of nodes to the isosceles triangle 
 $\Delta^I:=\{(k, l) \in \mathcal{P} :   k>0 \}$ and similarly for type II points : $\Delta^{II}:=\{(k, l) \in \mathcal{P}'  :   l>0 \}$. In    The union of $\Delta^I$ and $\Delta^{II}$ defines  
 a straight triangle  $\Delta$ (see figure 5 a)).
In summary : 
%\begin{figure}[!h]
%\fbox{ \includegraphics[width= 4 cm, height = 4 cm]{triangle.pdf}
%\includegraphics[width= 4 cm, height = 4 cm]{liss458pert.pdf},
%\includegraphics[width= 4 cm, height = 4 cm]{liss4-5-cheb.pdf}
% }
%\caption{ knot $L(2,3,5,0,.1,.2)$, projection $L(2,3,0,.1)$, and associated Chebyshev figure $C(2,3)$}
%\end{figure}
\begin{lem}\label{noeudlissajous} A generic Lissajoux planar curve  $L(n_1,n_2,\phi)$, with  $n_1$ and $n_2$ coprime and $\phi$ small, and    defined by 
\begin{equation}
L(n_1,n_2,\phi):\left( 
\begin{array}{c}
[0,1] \longrightarrow \mathbb{R}^2  \\
t\mapsto  \left(\cos 2\pi n_1 t, \cos\left[2\pi n_2\left( t +\phi\right)\right]\right)
\end{array}
\right)
\end{equation}
  has $ n_{L(n_1,n_2,\phi)} =2n_1n_2 -n_1 -n_2$ nodes parametrized by the integer points that lie in the interior of  the
  straight  triangle $\Delta$  : \\
$  J =   \Delta \cap \mathbb{N} \times \mathbb{N} :=\{ (k,l)  :  k>0, l>0, n_2 k  + n_1 l < 2n_1n_2\}$.\\
The  node  vector is by definition  $( t_j )_{j\in J}  \in [0,1[ ^{n_{L(n_1,n_2,\phi)}}$ where,\\ for each $(k,l) \in J$:
 \begin{equation}\label{typeI}
 t_{kl}   :=\left\{
 \begin{array}{cc}
  -\phi + \frac{1}{2}\left( -\frac{k}{n_1} + \frac{l}{n_2}     \right) &\quad  {\rm if } \quad   n_1 l > n_2 k \quad {\rm (type\  I)} \\
   \frac{1}{2}\left(  \frac{k}{n_1} - \frac{l}{n_2}    \right) & \quad {\rm if } \quad    n_1 l < n_2 k  \quad {\rm (type \ II)} 
 \end{array}
  \right.
\end{equation}
\end{lem}
 
\subsection{Prime knot  diagrams}\label{premier}

 \begin{figure}[!h]\label{billard}
\fbox{ \includegraphics[width= 4.35 cm, height = 4.35 cm]{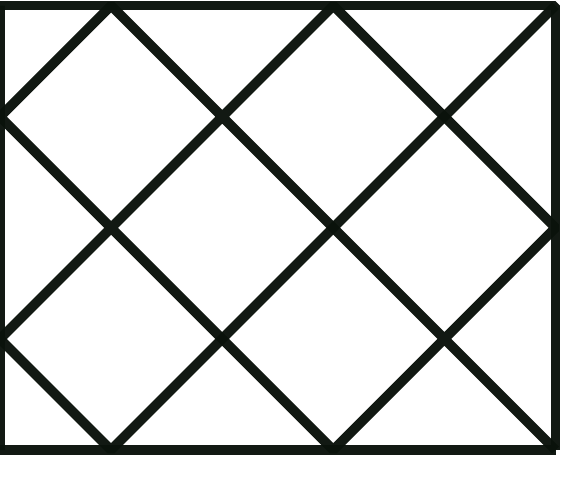}\hskip .5 in
\includegraphics[width= 4.3 cm, height = 4.4 cm]{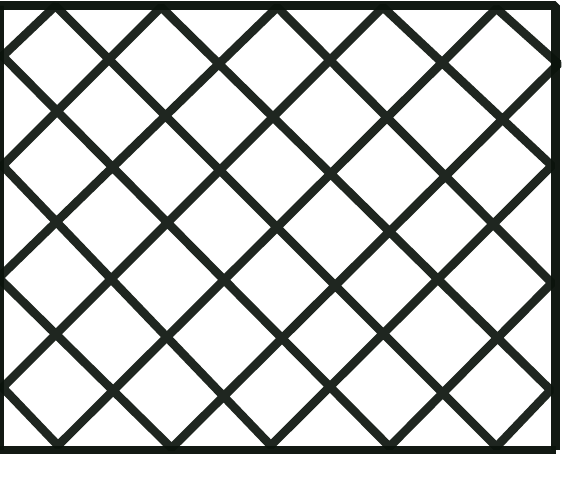}
\includegraphics[width= 5 cm, height = 5 cm]{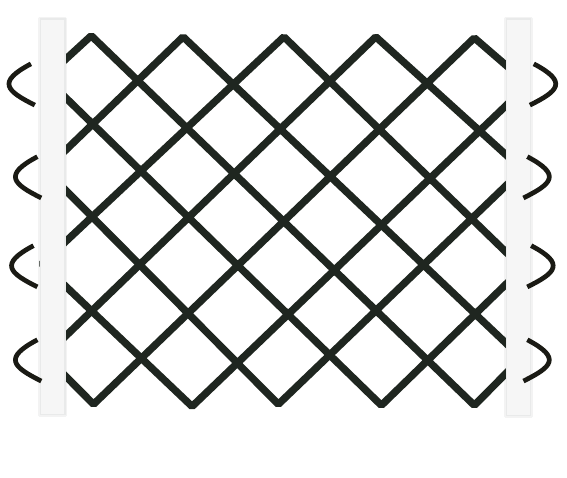}}
 \caption{   $C(4,5)$, billiard curve $L(4,5,.1)$ and   plat-closure  braid  of $L(4,5,.1)$ }
\end{figure}
We first  recall that a  checkerboard - such as $L(n_1,n_2)$-  is the plat-closure (cf. \cite{B}) of  a braid of $2n_1$ strands  defined by the following 
 group element of the braid group $B_{2n_1}$: 
 $$C( 2n_1,n_2-1) :=  \left( \sigma^*_2 \sigma^*_4 \cdots \sigma^*_{2n_1-2}\sigma^*_1 \sigma^*_3 \cdots \sigma^*_{2n_1-1}   \right)^{n_2-1}
 \sigma^*_2 \sigma^*_4 \cdots \sigma^*_{2n_1-2}$$
 (the powers of $\sigma_i$,  denoted  by$*$  are $\pm 1$; they  are irrelevant as far as the shadow is concerned). 

C. Lamm showed  in his thesis (cf \cite{La2} Theorem 2.3) that 
 any  knot   is presented by  a  checkerboard diagram  $C_k( 2n_1, 2n_1.p)$ for some $p$. 
This  diagram  has  a   Lissajous shadow of type $L(n_1,n_2=  2n_1p+1)$.  Hence we can  as well say that any knot
is presented by a Lissajous diagram.

 \begin{prop}\label{residu} Any knot  admits a Lissajous diagram with shadow $L(n_1,n_2) $ where $n_1$ and $n_2$ 
  are  odd primes  and where $  n_2  \equiv 1\ {\textrm mod\ } n_1 $.  
   \end{prop}
   \noindent {\it Proof :} 
 the only constraint on  $n_1$ is that it should be  greater or equal to the braid index of $K$. In particular choosing $n_1$ large enough,
we can suppose that  $n_1$  is  prime and odd. Moreover  from the proof of theorem 2.3 of \cite{La2}, we can also represent 
the same knot by diagrams with shadows $L(n_1,n_2=  2n_1(p+q)+1)$ for any nonnegative number $q$: 
$q$ represents the number of pure braids added to the original rosette braid. For each added piece 
and by a judicious choice of the crossing numbers 
 the new rosette braid  represents the same  knot.  By Dirichlet's prime number theorem and since $2n_1$ and $n_2$ 
are rel. prime,  there are numbers $q$ such that  $n_2$ is   prime.
   We will  thus restrict  our further  investigation of  the  Lissajous figure to the  case where $n_1$ and $n_2$ are odd primes.\\ 
  
%   \noindent {\bf Proof :} 
%   first statement follows from preceeding remarks. 
%   First equality of the congruence system  is a consequence of remarks in subsection \ref{symetriesdesnoeuds}.
  %  and  since choice of $n_3$ is free,  the existence of solutions of the  system of congruences
%    is deduced from the  chinese  remainder theorem. \hfill \qed\\
  %  There is an equivalent formulation for  toric shadows. We will only need :
   
    \subsection{Nodes coupling}\label{symetriesdesnoeuds}
     
  \begin{figure}[!h]\label{symmetry}
 \fbox{ \hskip .2 in  \includegraphics[width= 4 cm, height = 4 cm]{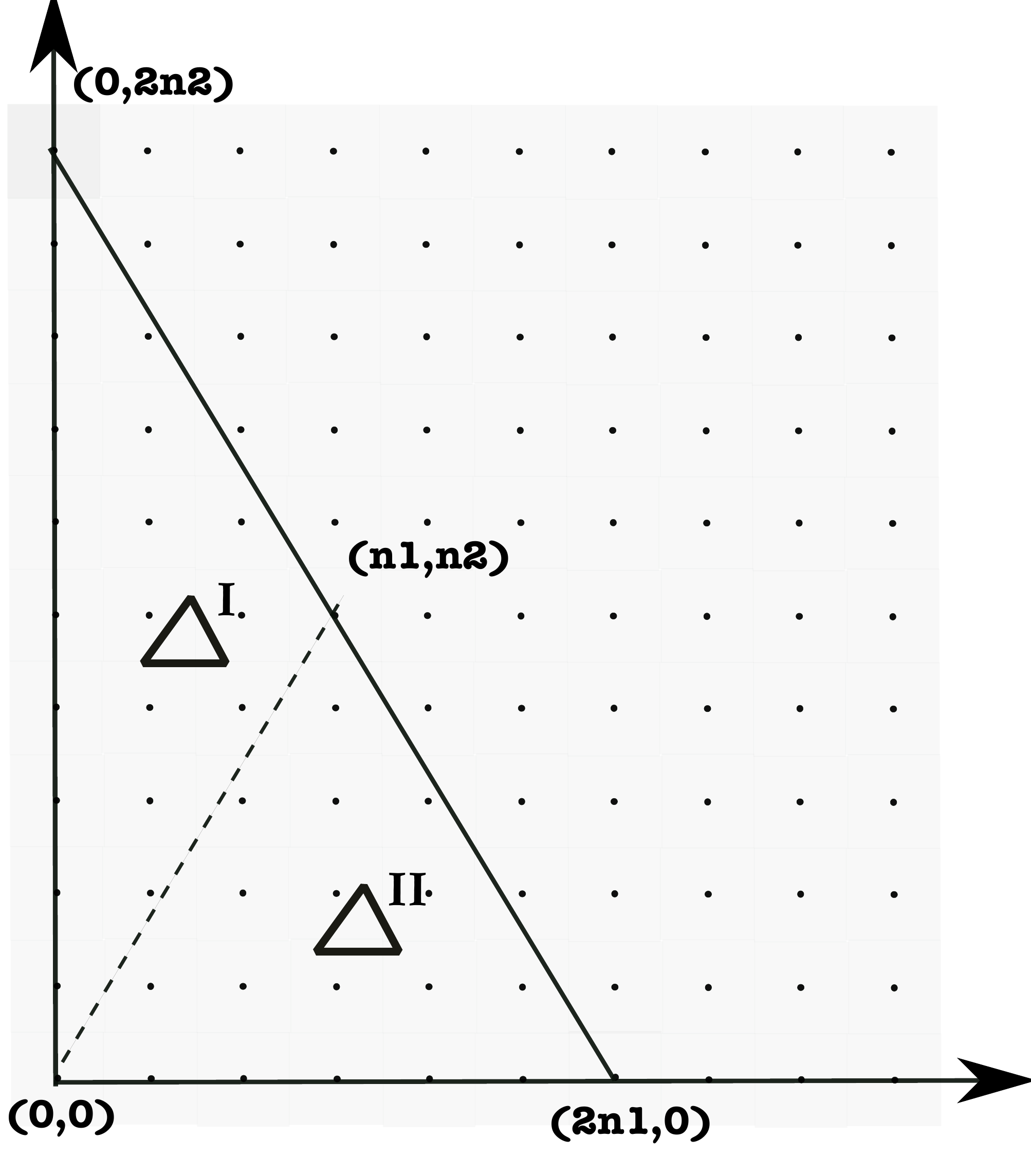}  \hskip .2 in 
\includegraphics[width= 4 cm, height = 4 cm]{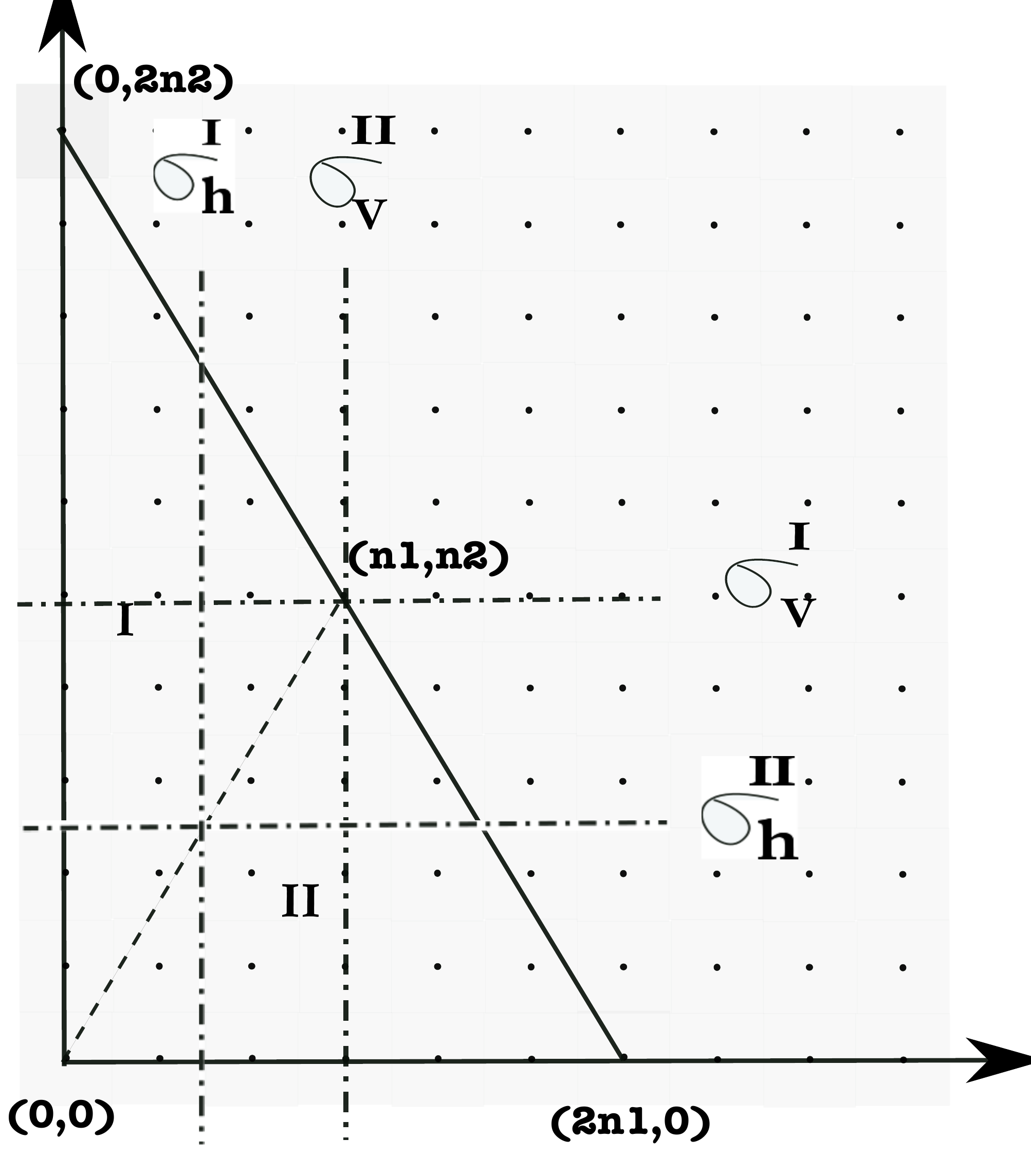} \hskip .2 in 
 }
\caption{Triangle $\Delta =\Delta^I \cup \Delta^{II}$;  hor. and  vert. symmetries}
\end{figure}

    %%%%%%%%%%%
    We describe a method to pair the nodes which will be crucial in section 5.\\
    
We will regroup   integer points of the triangle $\Delta$, described in Lemma \ref{noeudlissajous}, by pairs. The   $n_1n_2 -n_2$ nodes of type I (resp. the  $n_1n_2 -n_1$ nodes of type II)  are parametrized by the integer points of the
 isoceles subtriangle
 $\Delta^I : = \left( \left(0,0\right),\left(n_1,n_2\right), \left(0,2n_2\right)\right)\subset \Delta $\\
-resp.  
$\Delta^{II}: = \left( \left(0,0\right),\left(2n_1,0\right), \left(n_1,n_2\right)\right)  \subset \Delta $- (see fig 5) . 
     \begin{figure}[!h]\label{symmetry2}
 \fbox{  \hskip .2 in \includegraphics[width= 4 cm, height = 4 cm]{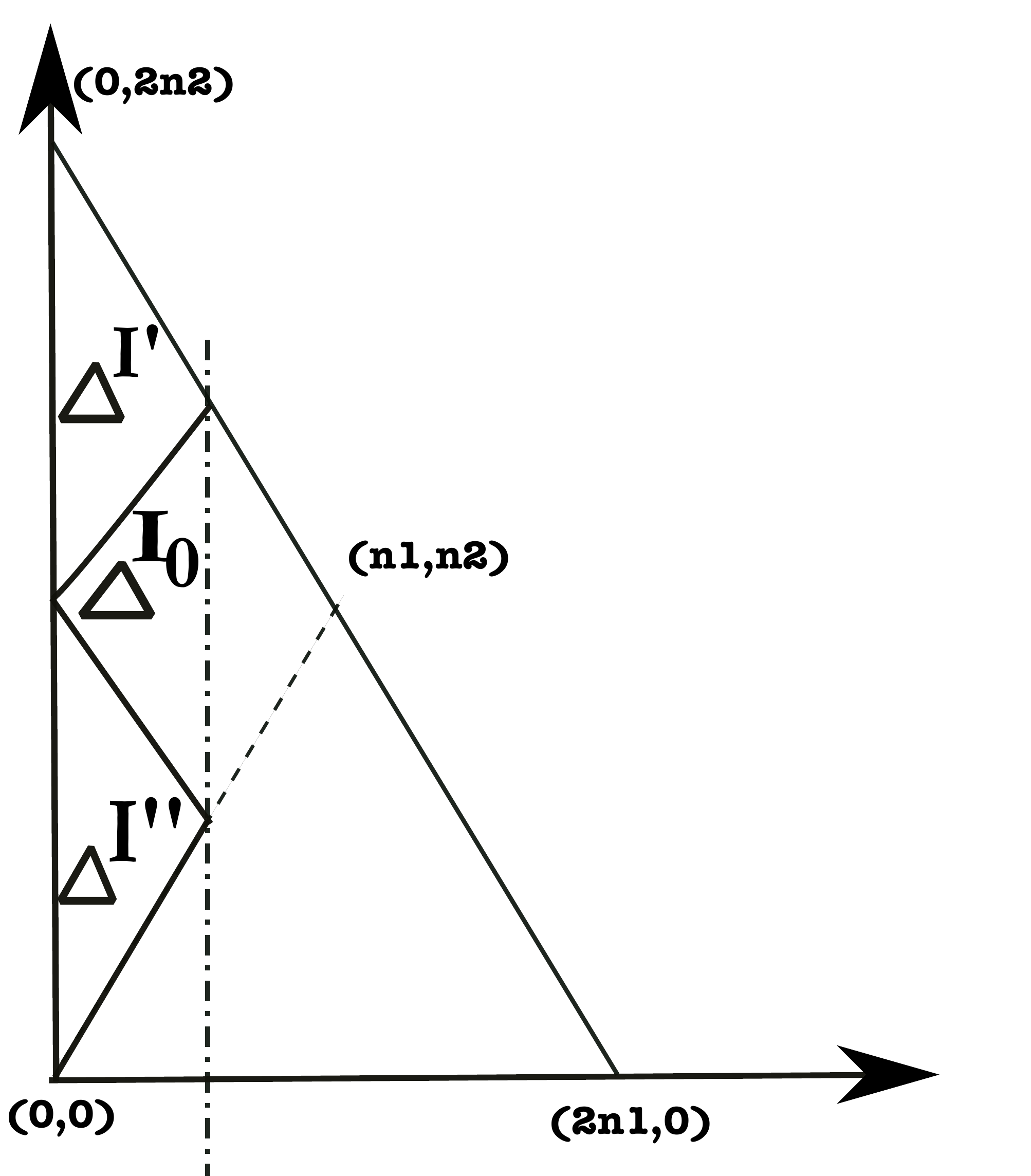}  \hskip .2 in 
\includegraphics[width= 4 cm, height = 4 cm]{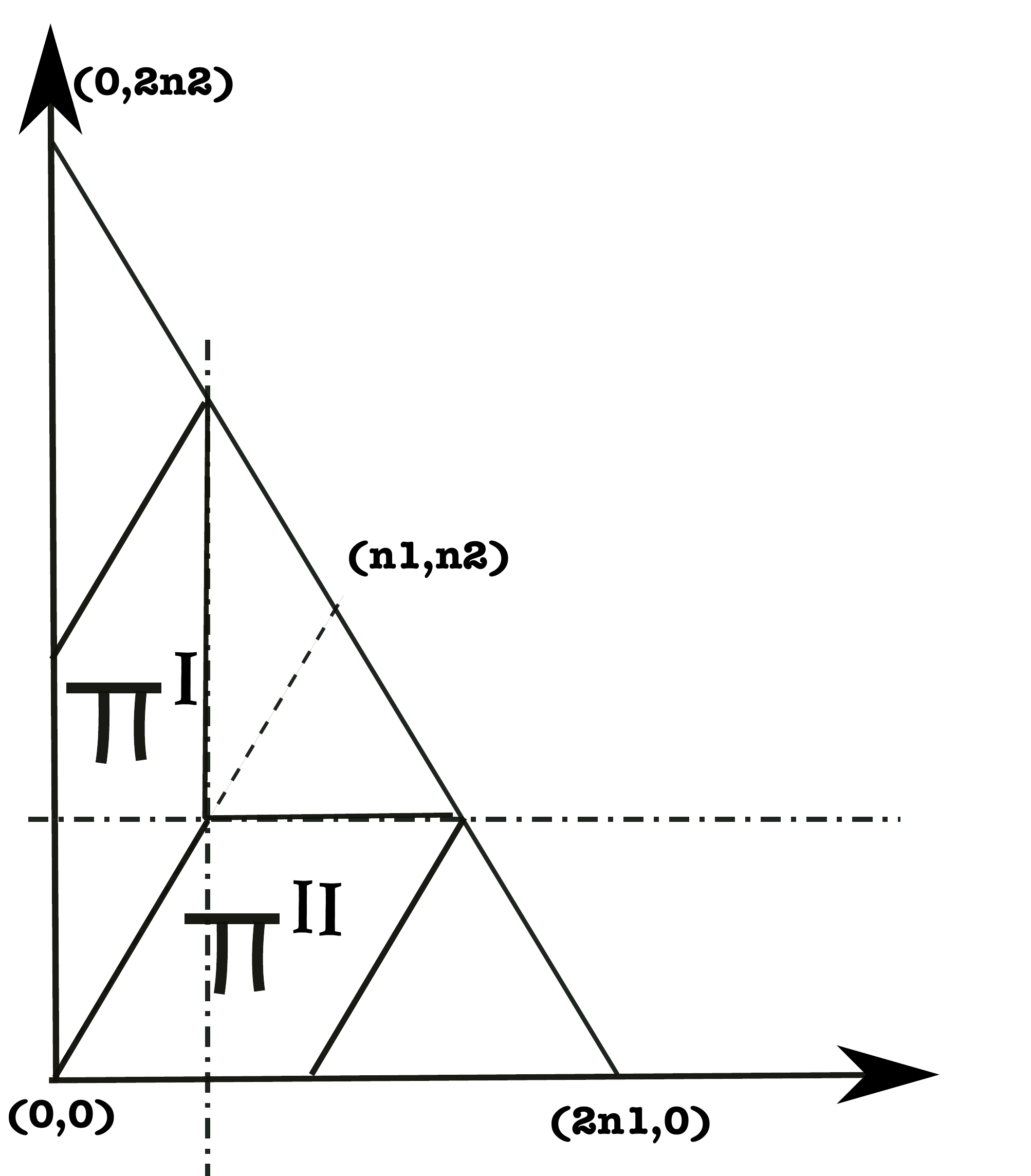} \hskip .2 in 
}
\caption{Subdomains of triangle $\Delta$}
\end{figure}
\begin{lem}\label{coupling} Each node of $\Delta^I$ (resp. $\Delta^{II}$)  is coupled to a different node of $\Delta^I$ (resp. $\Delta^{II}$).
The set of pairs extracted  from   $\Delta^I$ (resp. $\Delta^{II}$) are parametrized by integer points lying in  the parallelograms :\\
$\Pi^I := \left( \left(0,0\right),\left(\frac{n_1}{2}, \frac{n_2}{2}\right),
\left(\frac{n_1}{2}, n_2+\frac{n_2}{2}\right),\left( 0,n_2 \right)\right) $ \\  (resp.   
$\Pi^{II} := \left( \left(0,0\right),\left(\frac{n_1}{2},\frac{n_2}{2}\right),
\left(n_1+\frac{n_1}{2}, \frac{n_2}{2}\right),\left( n_1,0 \right)\right) $)(see figure 5).
\end{lem}
 \noindent {\it Proof :}
 \begin{enumerate}
 \item  We first regroup separately  integer points in  each isosceles triangle, namely  $\Delta^I $ and 
$\Delta^{II}$,  using resp.   the   horizontal  symmetry  $\sigma^I_h : (k,l) \longrightarrow (n_1-k, l)$ 
(resp. symmetry    $\sigma^{II}_h : (k,l) \longrightarrow (k,n_2-l)$).\\
We pair each  integer point above the axis of symmetry with  its  reflection image. Notice  that, since  $n_1$ and $n_2$ are odd,  no integer point lies 
on the axis of $\sigma^I_h$  (resp.   $\sigma^{II}_h$). Hence this process forms pairs of distinct points parametrized by 
integer points in the isosceles  triangle (see fig.6) \\
$\Delta^{I_0}:=\left(\left(0,n_2\right), \left(n_1/2,n_2/2\right),  \left(n_1/2,3n_2/2\right) \right)$.\\ 
(resp.  $\Delta^{II_0}:=\left(\left(n_1,0\right), \left(n_1/2,n_2/2\right),  \left(3n_1/2,n_2/2\right) \right)$).

  \item  The remaining points  of $\Delta^I$ that are not yet coupled, form two  isosceles subtriangles  $\Delta^{I'}, \Delta^{I''}$
  -resp  $\Delta^{II'}, \Delta^{II''}$ of $\Delta^{II}$-  (see figure 5). \\
  We then couple each point of $\Delta^{I'}$ with a point of $\Delta^{I''}$ (resp. 
 each point of $\Delta^{II'}$ with a point of $\Delta^{II''}$) via a translation $\tau^I$ of vector $(0,n_2)$
  (resp. via a translation $\tau^{II}$  of vector$(n_1,0)$ (see fig. 6).\\
 \end{enumerate}
 Alltogether the set of pairs of nodes   are parametrized by  the union of two parallelograms $\Pi^I :=\Delta^{I'} \cup 
   \Delta^{I_0}$
 and   $\Pi^{II} :=\Delta^{II'} \cup   \Delta^{II_0}$.
  \hfill \qed

 \subsection{Deformations of  planar Lissajous curves}
Let us consider a deformation of a Lissajous figure $\{ L_\epsilon(n_1,n_2,n_3,\phi,\psi  ) \}_{\epsilon\in [0,\epsilon_0[} $
   \begin{equation*}\label{plf}
L_\epsilon(n_1,n_2,n_3,\phi,\psi ):
\left(  
\begin{array}{ll}
[0,1]  &\longrightarrow \mathbb{R}^2 \\
t  &\mapsto \left( \cos\left(2\pi n_1 t\right ) ,   \cos\left(2\pi n_2\left( t +\phi\right)\right ) +\epsilon 
\cos\left(2\pi n_3\left(t+\phi +\psi\right) \right) \right) 
\end{array}
\right)  
\end{equation*}
 These are Fourier curves  of type $(1,2)$. The node positions   of the curve  are shifted as we can see in figure 6 and 
 the symmetry  of the initial Lissajous figure is  clearly broken.

\begin{figure}[!h]\label{deform}
\fbox{ \includegraphics[width= 4 cm, height = 4 cm]{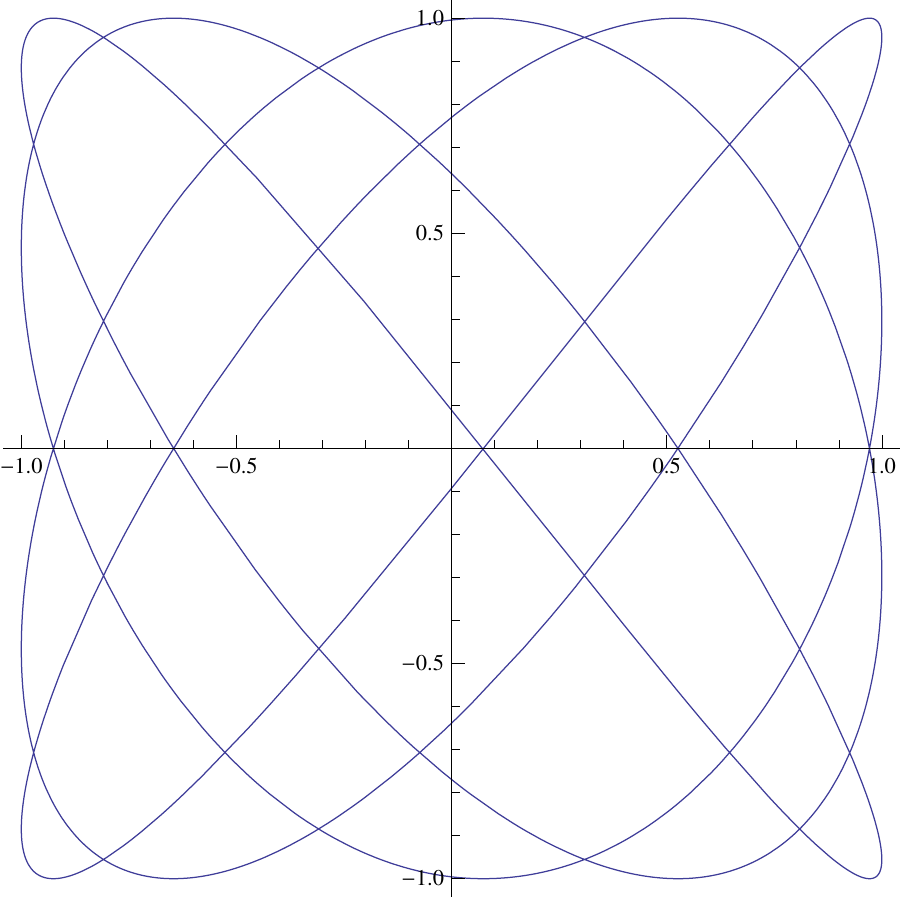}\hskip .5 in
\includegraphics[width= 4 cm, height = 4 cm]{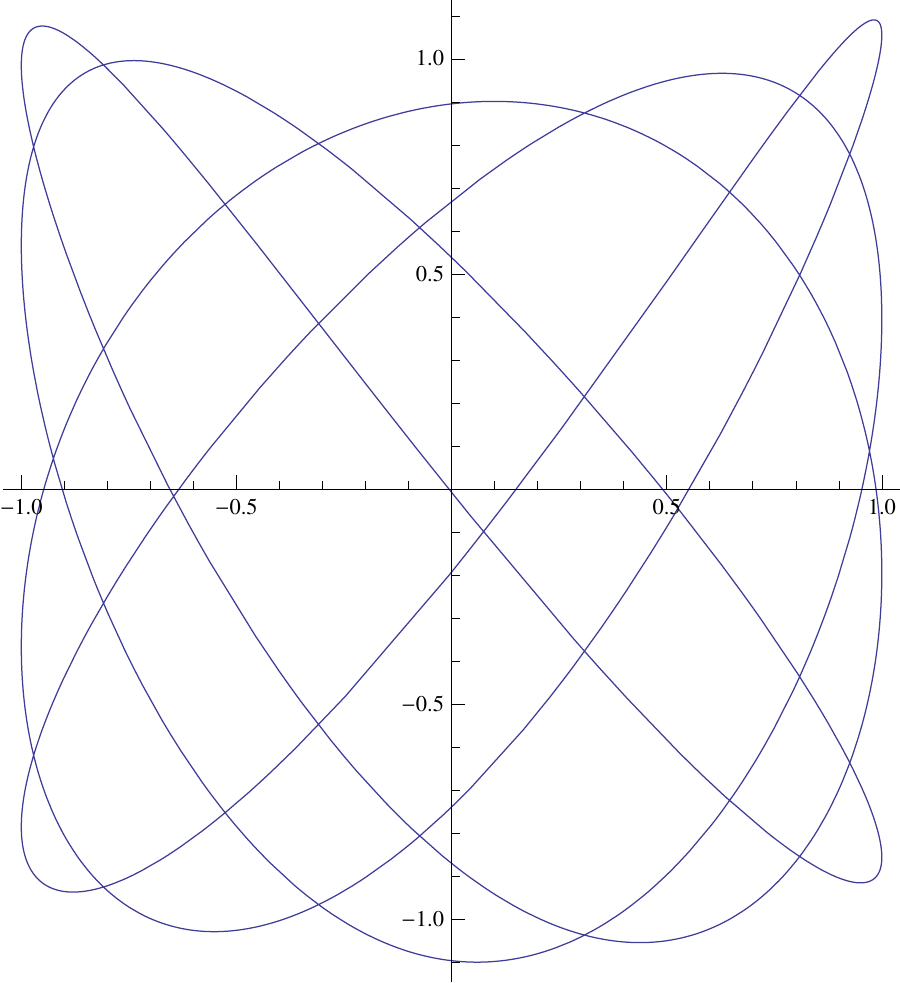}
 }
\caption{Lissajous shadow $L(4,5,.1)$ and deformation $L_{\epsilon}(4,5,7,.1,.2)$ }
\end{figure}
Let us  compute   the nodes parameters of these deformed Lissajous planar curves. 
\subsection{The nodal curve  of a    deformation of Lissajous curves}\label{nodallissajous}
Consider a deformation  $\{L_\epsilon(n_1,n_2,n_3,\phi,\psi )\}_{\epsilon \in [0,\epsilon_0]} $ 
  where  $n_1$ and $n_2$ are odd coprimes.\\
The coordinates of the nodal curve  $\eta   $ introduced  in section \ref{nodalcurve}
are ordered according to the lexicographic ordering of integer points of the 
triangle $\Delta$ described in Lemma  \ref{noeudlissajous}.We proceed as in this Lemma.\\
  
  Suppose $L_\epsilon(t) = L_\epsilon(s)$; then equality of first coordinates  yields 
 \begin{equation} \label{eq1bis}
   s = \sigma t + \frac{k}{n_1}, \quad k \in \mathbb{Z}, \quad \sigma = \pm 1.
   \end{equation} 
  Equality for  the second coordinates  yields : 
\begin{equation}\label{eq2}
\begin{array}{l}
\cos  \left(  2\pi n_2  \left( t +\phi \right) \right) - \cos  \left(  2\pi n_2 \left( s +\phi \right) \right)
=  \\ 
\epsilon \left(  \cos(2\pi n_3 (s+\phi +\psi) ) -  \cos(2\pi n_3 (t+\phi +\psi) )           \right)
\end{array}
\end{equation}
Or,  using  formula $\cos a -\cos b = -2\sin ( \frac{a+b}{2}) \sin (\frac{a-b}{2}),$
\begin{equation}\label{eq3zero}
\begin{array}{cc}
&\sin \left(  2\pi n_2  \left( \frac{t + s}{2}   +\phi \right) \right) \sin \left(  2\pi n_2   \frac{t - s}{2}   \right)  = \\
 & - \epsilon  \sin \left(  2\pi n_3  \left( \frac{t + s}{2} +\phi +\psi\right) \right) \sin \left(  2\pi n_3   \frac{t - s}{2} \right). 
\end{array}
\end{equation}

Plug  equality \eqref{eq1bis} into \eqref{eq3zero}; we get two cases according to the value of $\sigma=\pm 1$  :
\begin{enumerate}
\item{\bf Type I. } If  $\sigma =+1$ ; $ s =   t + \frac{k}{n_1} $ hence $\frac{s+t}{2}= t +\frac{k}{2n_1}$, $k\in\mathbb{Z}$  and  $t$ is given implicitely by the equation : 
 \begin{equation}\label{eq3}
\sin \left(  2\pi n_2  \left( t +\frac{k}{2n_1} +\phi \right) \right) 
a_k
=  \epsilon  \sin \left(  2\pi n_3  \left( t +\frac{k}{2n_1} +\phi +\psi\right) \right)
\end{equation}
where $ a_k = -  \frac{\sin   \left(  \pi\frac{n_2 k}{n_1}\right) }
{ \sin \left( \pi  \frac{n_3  k}{n_1}   \right)} $ (which is well-defined and not zero since $k<n_1$ and $n_2,n_3$ are prime
relatively to $n_1$).\\
Hence for any $\epsilon$,   solutions are  parametrized by  the  integer points  $(k,l)\in\Delta^I$ and given by 
\begin{equation}\label{typeIbis}
 t_{kl}(\epsilon) := 
   -\phi + \frac{1}{2}\left( -\frac{k}{n_1} + \frac{l}{n_2}     \right)   +  u_{kl}(\epsilon)
 \end{equation}
  where $u_{kl}$ is a function defined implicitely by  $f_{kl}\left(u_{kl}\left(\epsilon\right)\right) =\epsilon$ with : 
 $$f_{kl}( y )  = a_k(-1)^l\frac{\sin  2\pi n_2 y   }
 {\sin \left(  2\pi n_3  \left( y+  \frac{l}{2n_2} +\psi\right) \right)}$$

\item{\bf Type II.} If $\sigma =-1$ then   $ s =   -t + \frac{k}{n_1} $,  hence $\frac{t-s}{2} = t  -\frac{k}{2n_1}$ and 

 \begin{equation}\label{eq2}
 b_k\sin \left(  2\pi n_2  \left( t  -\frac{k}{2n_1} \right) \right) 
=  \epsilon    \sin \left(  2\pi n_3  \left( t  -\frac{k}{2n_1} \right) \right) 
\end{equation}
where $ b_k  = -  \frac{\sin \left(  2\pi n_2  \left(\frac{k}{2n_1}  +\phi \right) \right)}
{\sin \left(  2\pi n_3  \left( \frac{k}{2n_1} +\phi +\psi\right) \right) } $ 
($b_k$ is not zero and well-defined for small $\phi, \psi$  as for the $a_k$'s).
 \end{enumerate} 
 As in case I, solutions are given by 
\begin{equation}\label{typeIbis}
 t_{kl}(\epsilon) := 
   + \frac{1}{2}\left( \frac{k}{n_1} - \frac{l}{n_2}     \right)   +  u_{k,l}(\epsilon)
 \end{equation}
 with 
  $f_{kl}\left(u_{kl}\left(\epsilon\right)\right) =\epsilon$
  where 
 $$f_{kl}( y )  = b_k(-1)^l\frac{\sin  2\pi n_2 y  }
 {\sin \left(  2\pi n_3  \left( y+  \frac{l}{2n_2}\right) \right)}.$$
Notice that  $u_{kl}(0)=0$ and  that the $ t_{kl}(0)$ are the  solutions of type I and  II described in  Lemma \ref{noeudlissajous}.  
Let us summarize these facts and define the  nodal curve of the deformation:
 %The deformation of  Lissajous shadows define ipso facro 
%a deformation of Lissajous knot  These deformations are still Fourier knots of type $(1,1,2)$.  
%All knots are isotopic  to Fourier knots 
 %but Lissajous knots as well as Simple (minimal) knots offer strong symetries. 
\begin{lem}\label{deflis}
Consider a deformation   of planar Lissajous curves $L_\epsilon(n_1,n_2,n_3,\phi,\psi )$ with $n_1$ and $n_2$ odd   and consider a
Lissajous  deformation :
  \begin{equation*}\label{plf}
L_\epsilon(n_1,n_2,n_3,\phi,\psi ):
\left(  
\begin{array}{ll}
[0,1]  &\longrightarrow \mathbb{R}^2 \\
t  &\mapsto  \left(\cos 2\pi n_1 t , \cos\left[2\pi n_2\left( t +\phi\right)\right] +\epsilon 
\cos\left[2\pi n_3\left(t+\phi +\psi\right)  \right] \right)
\end{array}
\right)  
\end{equation*}
 $\epsilon \in [0,\epsilon_0[$,  where we suppose that $\epsilon_0$  and $\phi$ are   small and $n_1, n_2, n_3$ are relatively 
 prime.\\
 Each  planar curve  $L_\epsilon$  has $ n_{L(n_1,n_2,\phi)} =2n_1n_2 -n_1 -n_2$ nodes
  parametrized by the integer points   that lie in the interior of  the
straight   triangle $\Delta$ (defined in Lemma \ref{noeudlissajous})\\
The  nodal   curve  $ \eta (\epsilon) = ( t_{kl} )_{(k,l)\in \Delta}  \in [0,1[ ^{n_{L(n_1,n_2,\phi)}}$ 
is defined by  
 \begin{equation}\label{typeI}
\eta(\epsilon)_{kl} =  t_{kl}(\epsilon)  :=\left\{
 \begin{array}{cc}
  -\phi + \frac{1}{2}\left( -\frac{k}{n_1} + \frac{l}{n_2}     \right) + u_{kl}(\epsilon)   &\quad  {\rm if }\  (k,l)\in \Delta^I\  \\
   \frac{1}{2}\left(  \frac{k}{n_1} - \frac{l}{n_2}    \right) + u_{kl}(\epsilon)   & \quad  {\rm if }\  (k,l)\in \Delta^{II}.
 \end{array}
  \right.
\end{equation}
where the $kl$ are lexicographically ordered.\\
Functions $u_{kl}$  are the local  inverse functions of  functions $f_{kl}$, ie such that 
$  f_{kl}\circ u_{lk} = Id $.   \\
$\bullet$ If  $(k,l)\in \Delta^I$ then :
 $$    f_{kl}( y )  =\frac{\sin  2\pi n_2 y   }{ a^I_{kl}\sin \left[  2\pi n_3 
   \left( y+  \frac{l}{2n_2} +\psi\right) \right]}, \quad    a^I_{kl} : = (-1)^{l+1}  \frac{ \sin  \pi  \frac{n_3  k}{n_1}    }
   {\sin   \pi\frac{n_2 k}{n_1} } $$
 
\noindent  $\bullet$   If  $(k,l)\in \Delta^{II}$ then  :

 $$  f_{kl}( y )  =\frac{\sin  2\pi n_2 y  }
 { a^{II}_{kl}\sin \left[  2\pi n_3  \left( y+  \frac{l}{2n_2}\right) \right]},\qquad  a^{II}_{kl} :
  = (- 1)^{l+1} \frac{\sin \left[  2\pi n_3  \left( \frac{k}{2n_1} +\phi +\psi\right) \right] }{\sin \left[ 2\pi n_2  \left(\frac{k}{2n_1}  +\phi \right) \right]}.$$
 \end{lem}
%Notice that the $a_{kl}$ can be seen as polynomials in $\cos\psi$ and $\sin\psi$:
%$$a^{I}_{k,l} := a^{I}_{k,l,0} \cos \psi + b^{I}_{k,l,0} \sin\psi $$
%where 
%$ a^{I}_{k,l,0} := (-1)^{l+1}  \frac{ \sin \left( \pi  \frac{n_3  k}{n_1}   \right)\sin   \left(  \frac{2\pi n_3.l}{n_2}  \right)}
%   {\sin   \left(  \pi\frac{n_2 k}{n_1}\right) }, b^{I}_{k,l,0} := (-1)^{l+1}  \frac{ \sin \left( \pi  \frac{n_3  k}{n_1}   \right)\cos  \left(  \frac{2\pi n_3.l}{n_2}  \right)}
%   {\sin   \left(  \pi\frac{n_2 k}{n_1}\right) }$
%  And
%  $$a^{II}_{k,l} := a^{II}_{k,l,0} \cos \psi + b^{II}_{k,l,0} \sin\psi $$ 
%  where 
% $ a^{II}_{k,l,0} = (- 1)^{l+1} \frac{\cos \left(  2\pi n_3  \left( \frac{k}{2n_1} +\phi  \right) \right)\sin \left(    \frac{2\pi n_3 .l}{n_2} \right) }{\sin \left(  2\pi n_2  \left(\frac{k}{2n_1}  +\phi \right) \right)}, 
% b^{II}_{k,l,0} = (- 1)^{l+1} \frac{\sin \left(  2\pi n_3  \left( \frac{k}{2n_1} +\phi  \right) \right)\sin \left(    \frac{2\pi n_3 .l}{n_2} \right) }{\sin \left(  2\pi n_2  \left(\frac{k}{2n_1}  +\phi \right) \right)}$
%\end{rem}

\section{$\mathbb{Q}$-Linear dependence of  the nodal curve coordinates}
%Given any curve in $\gamma : [0,\epsilon_0] \longrightarrow \mathbb{R}^n$, we find a simple 
%geometrical property  that ensures that  
%a dense subset of $[0,\epsilon_0] $ the set of coordinates of 
%$\gamma (\epsilon) := \left(y_1 (\epsilon), \cdots, y_n (\epsilon),1   \right) $ are rationally linear dependent.
  \subsection{Infinitesimal deformation and  $\mathbb{Q}$-linear independence}
  We will give a infinitesimal  condition on the curve to ensure 
 the rational linear independence. But first let us prove 

\begin{lem}\label{courbe}Let $\gamma : I = [0,\epsilon_0] \longrightarrow [0,1]^n$ be a real analytic  skew curve. 
Then the set of $\epsilon \in  [0,\epsilon_0] $ such that   the    numbers 
\begin{equation}\label{nombres}
 \gamma_1(\epsilon), \cdots,   \gamma_n(\epsilon), 1   
 \end{equation} 
are rationally  linear independent,  is dense in 
$[0,\epsilon_0] $.
\end{lem}
\noindent{\it Proof :} 
we exhibit a sequence of $\epsilon$'s that converges to zero and such that the  numbers \eqref{nombres}
are rationally linear independent-the same proof  will work for any other points $\neq 0$.\\
 If such a sequence does not exist,  there is a  nonempty neighborhood $ I'$ of $0$ in $I$  such that 
the coordinates of $\gamma(\epsilon)$ and $1$ are rationally  linear dependent; then 
$$\sum_{i=1}^n \lambda_i(\epsilon) \gamma_i (\epsilon ) + \lambda_{n+1}(\epsilon) = 0.$$
Since $I'$ is an uncountable set, the map  $ \lambda: I'  \longrightarrow \mathbb{Q}^{n+1}$  takes    the  same vector 
 value  $(\lambda_i)_{i=1,\cdots,n}$
 for infinitely many $\epsilon_k$ with 
$\underset{k\mapsto+\infty}{\lim} \epsilon_k = 0$. Hence for infinitely many values of $\epsilon_k$ converging to $0$ 
we have 
\begin{equation}\label{relation}
\sum_{i=1}^n \lambda_i\gamma_i (\epsilon_k ) + \lambda_{n+1}= 0.
\end{equation}
Since $\gamma$ is analytic we deduce that this relation is true for all $\epsilon$; hence $\gamma$ belongs to 
a hyperplane with rational coefficients, which contradicts the skewness of curve $\gamma$.
 \hfill \qed\\
Similarily we prove 
 \begin{lem}\label{dephase} If $t_1,\cdots, t_n, 1$ are rationally lin. independent, then 
 $\{t_1+\phi ,\cdots, t_n+\phi, 1\}$  are rationally linear  independent for a dense subset of $\mathbb{R}.$ 
 \end{lem}
Taking successive derivatives of  equation   \eqref{relation} in Lemma \ref{courbe}
we deduce the  following infinitesimal   criterion for $\mathbb{Q}$-linear independence :\\
\begin{cor}\label{wronskien} Let $\gamma :[0,\epsilon_0] \longrightarrow [0,1]^n$ be an analytic curve   
 such that a   Wronskian is not zero :  
 \begin{equation}
   \left|
   \begin{array}{cccc}
   \gamma_1^{(k_1)}(0) &  \gamma_1^{(k_2)}(0)  & \hdots  &   \gamma_1^{(k_n)}(0)    \\
      \gamma_2^{(k_1)}(0) &  \gamma_2^{(k_2)}(0)  & \hdots  &   \gamma_2^{(k_n)}(0)    \\
   \vdots  & \vdots & \hdots  & \vdots     \\
   \gamma_n^{(k_1)}(0)      &  \gamma_n^{(k_2)}(0)  &  \hdots  & \gamma_n^{(k_n)}(0)    
   \end{array}
    \right| \not = 0
   \end{equation}
where $k_1<k_2<\cdots, k_n$; then  the set of real numbers  $\{1, \gamma_j (\epsilon) , \ j =1,\cdots, n \}$ are linear independent 
 over $\mathbb{Q}$ for a dense subset of $[0,\epsilon_0]$.
\end{cor}

  \subsection{The Lissajous nodal curve is skewed}\label{a practical case}
 We will show in this section that 
  \begin{prop}\label{proposition1} 
  Assume that 
   $n_1,n_2$ are  odd primes  and $n_3$ is  relatively prime to $n_1,n_2$ and  satisfy :
   \begin{equation}
  \left\{
  \begin{array}{cc}
    n_3  &\equiv 2 \quad {\textrm mod\ } n_1\\
      n_3  &\equiv 2 \quad {\textrm mod\ } n_2.
    \end{array}
    \right.
    \end{equation}
   Then there are positive 
  real numbers $\epsilon_0, \phi_0,\psi_0$ such that   for any 
  $|\epsilon|\leq \epsilon_0$, $|\phi  |\leq \phi_0, |\psi-\frac{\pi}{4}| \leq \psi_0 $,  the nodal curve $\eta(\epsilon)$ of the Lissajous deformation   $L_\epsilon(n_1,n_2,n_3 ,\phi,\psi )$  is skewed.
   \end{prop}
  In Subsection \ref{sub1} we give an explicit expression of the nodal coordinates of a  Lissajous deformation;
  we  then show, in  Subsection \ref{sub2},  that the Wronskien $D$  of the nodal curve is a polynomial in $r =\frac{n_2}{n_3}$;
  in Subsection \ref{sub3} we compute the coefficient $D_0$ of the lowest order monomial.
  \subsection{General form of the nodal curve coordinates}\label{sub1} 
 Consider   a curve $Y:$    
\begin{equation}
  Y : \left( 
  \begin{array}{cc}
   I:= [0,\epsilon_0] &\longrightarrow \mathbb{R}^n\\ 
     \epsilon &\mapsto    \left(  y_j\right)_{j \in J}
    \end{array}
    \right)
    \end{equation}
    with initial position $Y(0) = \left( y_j\left(0\right)\right)_{j=1,\cdots, n}$ and let us define curve $U$ such that 
    $Y = Y(0) + U.$
Then  :
  $$  y_j(\epsilon)  = y_j(0) + u_j(\epsilon) \ {\rm and} \ u_j(0)=0 \ \quad
   \forall \epsilon\in I,   j =1, \cdots, n.$$
  For each $ j=1,\cdots, n$, the $j$-th coordinate function $u_j $ of $U$ is defined as  the inverse function of a given function   $f_j$ : 
  $$f_j (u_j(\epsilon)) =\epsilon \quad  \forall\epsilon \in I.$$
   Though  each coordinate function $u_{j}(\epsilon)$- which we will denote  indistinctively by $u(\epsilon)$-  is    defined  implicitely,  the Lagrange inversion formula (cf. for instance \cite{S}) yields an explicit    power series expansion in terms of $\epsilon $  in the neigborhood of any 
   $0$ :
        $$u (\epsilon) :=  \sum_{n\geq 1} u_{n}(0) \frac{\epsilon^n}{n!}$$
  Namely :
  \begin{equation}\label{lagrange}
  u(\epsilon) =   \sum_{n=1}^\infty \frac{d^{(n-1)}}{dt^{n-1}}
  \left(\frac{t}{f(t)}\right)^n\biggr|_{t=0} \frac{\epsilon^n}{n!}.
    \end{equation}
 
 \subsection{Series expansion of   nodal curves}\label{sub12} 
Each one of the coordinates of  the nodal   curve of   a  Lissajous    deformation
 is defined by the inverse  of a real    function  of the following type :  
\begin{equation}\label{fonctionf}
f (u) = \frac{\sin u}{a\sin (r u+\psi)}\qquad  (  r\in\mathbb{Q}, a\in \mathbb{R}^*,  \psi \in \mathbb{R}). 
\end{equation}
 A simple computation  using  Lagrange's inverse formula \eqref{lagrange} (cf. for example \cite{S})  gives  the coefficients $u_{n}(0)$
  of the expansion  around zero of $u$ such that $f\circ u = Id$:
        $$u_{1}(0) = a\sin\psi ,  u_{2}(0)= 
        a^2r\sin\psi \cos\psi,, u_{3}(0)= a^3\sin\psi\left(6r^2 +(1-9r^2)\sin^2\psi\right)\cdots .$$
      
     \begin{lem}\label{calcul} Let $u$ be the inverse function of $f$  defined by equation  \eqref{fonctionf}. \\
      For each $n\geq 1$ the n-th coefficient of the power series of $u$ around zero 
        :  $\ u_{n}$   equals $a^n P_n ( \sin\psi, \cos\psi, r) $, where 
      $P_n$ is  a polynomial  of degree $n$   w.r.t.   variables $\cos\psi$ and $\sin\psi$. 
     $P_n$ is  also a  polynomial in  $r$ of degree $n-1$ whose lowest order term is:
      \begin{equation}\label{loworder}
 \bar u_{n}(0) :=\left\{
 \begin{array}{ll}
 a^n\left(\frac{y}{\sin u}\right)^n\biggr|_{u=0}^{(n-1)} \sin^n \psi &:=
   a^n c_n \sin^n \psi\quad \hfil ( n\  {\rm odd}), \\
   r.n a^n(\frac{y}{\sin u})^{n}\biggr|_{u=0}^{(n-2)} \cos \psi \sin^{n-1} \psi  &:=   a^n c_n \cos \psi \sin^{n-1} r\psi \quad  {\rm  (n  \ even)}
 \end{array}
  \right.
 \end{equation}
        \end{lem}
        {\it Proof :}    Lagrange's  inversion formula \eqref{lagrange} applied to function \eqref{fonctionf} and evaluated at zero, yields : 
                      $$ u_{n}(0) = a^n \left(  \left(\frac{y}{\sin y}   \right)^n \sin^n\left(r y +\psi     \right)     \right)\biggr|_{y=0}^{(n-1)}$$
%
 %           $$ u_{,n}(\epsilon_1) = a^n \frac{d^{(n-1)}}{dt^{n-1}} \biggr|_{t=\mu_1}  \left(\frac{t-\mu_1}
 %               {\frac{\sin t}{\sin\left(r t +\varphi \right) } - \frac{\sin \mu_1}{\sin\left(r \mu_1 +\varphi \right) }}   \right)^n $$
 %   Since 
 %   $$\frac{t-\mu_1}{\frac{\sin t}{\sin\left(r t +\varphi \right) } - \frac{\sin \mu_1}{\sin\left(r \mu_1 +\varphi \right) }}
%        =  \frac{(\cos r\mu_1\cos rt\sin\varphi  )(t-\mu_1)}{\cos r\mu_1\sin t - \cos r t \sin \mu_1\right)} }+o(t,\mu_1)$$
   \noindent Leibniz formula yields :
         $$ u_{n}(0) = a^n
         \sum_{k=0}^{n-1}  \binom{n-1}{k}\left[\left(\frac{y}{\sin\ y}   \right)^n\right]\biggr|_{y=0}^{(k)}
          \sin^n\left(r y +\psi     \right)\biggr|_{y=0}^{(n-1-k )}$$
          Since $\frac{y}{\sin y}$ is an even function, 
          $$ u_{n}(0) = a^n
         \sum_{k=0}^{\left[  \frac{n-1}{2}      \right]}  \binom{n-1}{2k}\left[\left(\frac{y}{\sin\ y}   \right)^n\right]\biggr|_{y=0}^{(2k)}
          \sin^n\left(r y +\psi     \right)\biggr|_{y=0}^{(n-1-2k )}$$
           $y_n$ is a polynomial in  the two variables $\cos\phi , \sin\phi$ and coefficients  are of the form $p_k(r)$- an  integer polynomial in $r$.  \\
  %            The  monomial coefficient  of highest order wrt r , ie $r^{n-1}$ is :
   
%\begin{equation}\label{formule1}
 % \left(\sin^n x \right) \biggr|^{(n-1)}_{\varphi}   = 
% \frac{1}{2^{n-1}} \sum_{k=0}^{\left[\frac{n}{2}\right]^*}  \binom{n}{k} (-1)^{k}(n-2k)^{n-1}\sin (n-2k)\varphi
% \end{equation}
% Hence 
%\begin{equation}\label{formule1}
%  Q_{n-1}(\cos\varphi)   =  
% \frac{1}{2^{n-1}} \sum_{k=0}^{\left[\frac{n}{2}\right]^*}  \binom{n}{k} (-1)^{k}(n-2k)^{n-1}
% U_{n-2k-1}(\cos \varphi)
% \end{equation}
%where $U_{n-1}(x)$ denotes  the Chebyshev polynomial of second order  defined by  
% $U_{n-1}(\cos\varphi) :=\frac{\sin n\varphi }{\sin \varphi }  $.\\
%%\begin{equation}\label{loworder}
 %\bar y_{,n} :=\left\{
% \begin{array}{ll}
% a^n\left(\frac{y}{\sin y}\right)^n\biggr|_{0}^{(n-1)} \sin^n \psi &=
%   a^n c_n \sin^n \psi\quad  {\rm  if\  n\ is \ odd}, \\
%   n(n-1)a^n(\frac{y}{\sin y})^n\biggr|_{0}^{(n-2)} \cos \psi \sin^{n-1} \psi  &=   a^n c_n \cos \psi \sin^{n-1} \psi \quad  {\rm  if\  n\ is \ even}
% \end{array}
%  \right.
 %\end{equation}
 The  monomial of lowest  order w.r.t. variable  $r$, is  of degree zero or one, according to the parity of $n$ and the corresponding coefficients are given by equation \eqref{loworder}.

\subsection{Series expansion of the  nodal curve's Wronskian  }\label{sub2}

The   coordinates of considered  nodal curves may be  ordered as follows.\\
  Consider  first  a set of roots  functions $u_m$  of  $ e $ equations in $u$
     $$ h_m\circ u_m = Id$$ 
     where 
    $$h_m ( u) := \frac{\sin u}{a_m\sin (r u+ \phi_m)}, \ m = 1,\cdots  e.$$
  We now specialize to integer  parameters  $(k,l)$ describing nodal curves.
   For each $k= 1,\cdots  e$,  consider  $n_k$   functional equations : 
    $$ f_{kl}\circ u_{kl}= Id \quad l = 1,\cdots , n_k$$ 
    where 
    \begin{equation}\label{definitiondef}
    f_{kl} ( u) :=     \frac{\sin (u )}{a_{kl}\sin (r u+ \psi_{kl})}.
    \end{equation}
    In our cases, inverse  functions $u_{kl}$  are indexed by a subset of 
      $J = \{ (k,l) :  1\leq k \leq e, \ 1\leq l \leq n_k \leq q  \}$; we then order $J$  according to 
  lexicographic order and construct the nodal  curve $(\eta^j(\epsilon))_{j\in J}$ accordingly.
  \begin{ex}For Lissajous curves $e=2n_1, q=2n_2$, $m= 2n_1n_2 -n_1-n_2 := 2p$,  $J= \Delta.$
    \end{ex}
Let us now   consider a  nodal curve of a deformation whose nodes are parametrized by set $J$   
($\#J = m$).
    \begin{equation}
    \eta : 
    \left( 
    \begin{array}{ll}
     [0,\epsilon_0] &\longrightarrow  \mathbb{R}^{ m } \\
      \epsilon  &\mapsto    \left(u^j\left(\epsilon\right)\right)_{j=1,\cdots,m} 
    \end{array}
    \right),
    \end{equation}
    and compute  its Wronskien $D$ at $t=0$:
        \begin{equation}
     D :=  | \eta'(0) , \cdots , \eta^{(m)}(0) | = \left|
   \begin{array}{cccc}
   u_{11}^{(1)} &  u_{11}^{(2)}  &   \hdots \hdots  \hdots    &   u_{11}^{(m)}    \\
     \vdots  & \vdots & \vdots\hdots \hdots  \hdots \vdots  & \vdots     \\
      u_{1n_1}^{(1)} &  u_{1n_1}^{(2)}  & \hdots \hdots  \hdots  &   u_{1n_1}^{(m)}    \\
        u_{21}^{(1)} &  u_{21}^{(2)}  & \hdots \hdots  \hdots   &   u_{21}^{(m)}    \\
   \vdots  & \vdots & \vdots\hdots \hdots  \hdots\vdots   & \vdots     \\
    u_{2n_2}^{(1)} &  u_{2n_2}^{(2)}  & \hdots \hdots  \hdots    &   u_{2n_2}^{(m)}    \\
    \vdots  & \vdots & \vdots \hdots \hdots \hdots \vdots   & \vdots     \\
     u_{e1}^{(1)} &  u_{e1}^{(2)}  &\hdots \hdots  \hdots    &   u_{e1}^{(m)}    \\
      \vdots  & \vdots & \vdots \hdots \hdots \hdots \vdots   & \vdots     \\
   u_{en_e}^{(1)} &  u_{en_e}^{(2)}  & \hdots \hdots  \hdots    &   u_{en_e}^{(m)}    \\  
   \end{array}
    \right|(0)
    \end{equation}
 \subsection{The Wronskian is not zero }
Our ultimate goal is to show that the Wronskian is not zero. However a direct computation is hopeless. We will overcome this difficulty 
  using two main  tricks. 
  First,   from Lemma \ref{calcul},  we can expand the determinant $D$  as  a polynomial in $r$; using  expansions \eqref{loworder}
    and  indexation given by \eqref{definitiondef}, we  can compute the coefficient $D_0$  which denotes 
   the   coefficient   of the monomial of $D(r)$ of lowest order w.r.t. $r$ (which is of degree $\frac{p(p+1)}{2}$).
   We will show that the non nullity of $D_0$ implies the non nullity of $D$ for almost all $r$.\\
   The second trick is that  $D_0$ can be computed and reduces to factors of a Vandermonde determinant when rows- that each corresponds to a node-  are grouped by pairs  just as was explained in Subsection \ref{symetriesdesnoeuds}.  \\
   First let us compute :   
       \begin{equation*}
  D_0 = \prod_{n=1}^m  c_n \left|
   \begin{array}{ccccc}
a_{11}\sin\psi_{11} &  a_{11}^2 \sin\psi_{11}\cos\psi_{11}  &  a_{11}^3 \sin^3\psi_{11} & a_{11}^4 \sin^3\psi_{11}\cos\psi_{11}      &       \hdots         \\
     \vdots  & \vdots & \vdots &\hdots  \hdots  &\vdots      \\
a_{1n_1}\sin\psi_{1n_1} &  a_{1n_1}^2 \sin\psi_{1n_1}\cos\psi_{1n_1}  &  a_{1n_1}^3 \sin^3\psi_{1n_1}        &  a_{1n_1}^4\sin^3\psi_{1n_1}\cos\psi_{1n_1} 
&  \hdots       \\
a_{21}\sin\psi_{21} &  a_{21}^2 \sin\psi_{21}\cos\psi_{21}  &  a_{21}^3 \sin^3\psi_{21} & a_{21}^4\sin^3\psi_{21}\cos\psi_{21}      &       \hdots        \\
    \vdots  & \vdots & \vdots   &\hdots \hdots \vdots       \\
a_{2n_2}\sin\psi_{2n_2} &  a_{2n_2}^2 \sin\psi_{2n_2}\cos\psi_{2n_2}  &  a_{2n_2}^3 \sin^3\psi_{2n_2} & a_{2n_2}^4 \sin^3\psi_{2n_2}\cos\psi_{2n_2}      &       \hdots        \\
   \vdots  & \vdots & \vdots  &\hdots \hdots \vdots       \\
      \vdots  & \vdots & \vdots  & \hdots \hdots \vdots        \\
a_{e1}\sin\psi_{e1} &  a_{e1}^2 \sin\psi_{e1}\cos\psi_{e1}  &  a_{e1}^3 \sin^3\psi_{e1}        &  \hdots    \hdots    &    \    \\
     \vdots  & \vdots & \vdots\hdots &\hdots  \hdots \vdots  & \vdots     \\
a_{en_e}\sin\psi_{en_k} &  a_{en_e}^2 \sin\psi_{en_k}\cos\psi_{en_k}  &  a_{en_e}^3 \sin^3\psi_{en_k}        &  \hdots    \hdots    &        \\
    \end{array}
   \right|
   \end{equation*}
   As the size   $m$  of the determinant $D_0$ is   even (m= 2p)  the   last column  consists  of terms of type $a^{2p}_{ij}\sin^{2p-1}\psi_{ij}\cos\psi_{ij}$.\\
   Let us introduce notations   $\alpha_{kl} := a_{kl}\sin\psi_{kl}$ and $\beta_{kl} := a_{kl}\cos\psi_{kl}$
   Then $D_0$ equals :
   \begin{equation}\label{dzero}
    D_0 =\pm \left(\prod_{n=1}^m  c_n  \right)\left(\prod_{(k,l)\in J}\alpha_{kl}\right) D_1
    \end{equation}
    where
     \begin{equation}\label{dzerun}
    D_1 := \biggl| 1, \alpha_{kl}^2, \cdots, \alpha_{kl}^{2p-2}, 
     \beta_{kl}, \beta_{kl}\alpha_{kl}^2, \cdots,  \beta_{kl} \alpha_{kl}^{2p-2}  \biggr|.
   \end{equation}
 
  We will compute $D_1$ for some specific values of $\phi,\psi$.
  \begin{lem}\label{dezerosept} Let us assume  $n_2$ is odd   and $n_3$ is even. Suppose also that $2\pi n_3\psi = \pi/4$ and $\phi=0$ ; then 
  $$D_1 = \pm  2^{m}\left(\prod_{(k,l)\in J}\beta_{kl} \right)
  \zeta^2(  \alpha^2_{kl} )_{(k,l)\in \Pi^{I} \cup \Pi^{II}}  $$
  where $\zeta (x_i)_{i\in I}:= \underset{i,j\in I, i<j}{\prod_{}}( x_i-x_j)$ is the difference product of the $(x_i)_{i\in I}$.
  \end{lem}
  \noindent {\bf Proof :}  
 For the sake of clarity we will denote  differently  $\alpha_{k,l}$ and $\beta_{k,l}$  according to whether 
 $(k,l) \in \Delta^I$ or
 $(k,l) \in \Delta^{II}$
 using coefficients  $a^I_{kl}$ and  $ a^{II}_{kl}$ as defined in Lemma \ref{deflis}   :
  
\begin{equation}\label{eq5} 
\left\{
 \begin{array}{ll}
 \alpha^I_{k,l}  =&(-1)^{l+1}
   \frac{ \sin  \pi  \frac{n_3  k}{n_1}    }{\sin   \pi\frac{n_2 k}{n_1} }   \sin \left[  2\pi n_3 \left(   \frac{l}{2n_2} +\psi\right) \right]\\
\beta^I_{k,l} =& (-1)^{l+1}  \frac{ \sin   \pi  \frac{n_3  k}{n_1}  }
   {\sin  \pi\frac{n_2 k}{n_1} }\cos \left[  2\pi n_3 \left(  \frac{l}{2n_2} +\psi\right) \right]     \\
 \alpha^{II}_{k,l}  =&(- 1)^{l+1} \frac{ \sin   \frac{\pi n_3 l}{n_2} }{\sin \left[ 2\pi n_2  \left(\frac{k}{2n_1}  +\phi \right) \right]}   
\sin \left[  2\pi n_3  \left( \frac{k}{2n_1} +\phi +\psi\right) \right]\\
  \beta^{II}_{k,l} =& (- 1)^{l+1} \frac{ \cos     \frac{\pi n_3 l}{n_2} }{\sin \left[  2\pi n_2  \left(\frac{k}{2n_1}  +\phi \right) \right]}   \sin \left[  2\pi n_3  \left( \frac{k}{2n_1} +\phi +\psi\right) \right] 
\end{array}
   \right.
  \end{equation}
 
 %We keep highest order term in $\cos\psi$; then coefficient of type $I$ and $II$ are identical up to variable $\phi$: 
 
%\begin{equation} 
%\left\{
%\begin{array}{ll}
% \alpha^{II}_{k,l}  &=(-1)^{l+1}
 %  \frac{ \sin  \pi  \frac{n_3  k}{n_1} \sin  2\pi n_3\left( \frac{l}{n_2}+\phi\right)    }{\sin \left[  2\pi n_2  \left(\frac{k}{2n_1}  +\phi \right) \right]}   \cos   2\pi n_3\psi
 %  =    \bar a_{k,l}(\phi)  \cos    2\pi n_3 \psi   \\
 %  \alpha^{I}_{k,l}  &= \bar a_{k,l}(0)  \cos   2\pi n_3 \psi     \\
 %    \beta^{II}_{k,l} &= (- 1)^{l+1} \frac{ \cos \frac{2\pi n_3 l}{n_2} 
%      \sin   2\pi n_3 \left( \frac{ k}{2n_1} +\phi  \right)}{\sin \left[  2\pi n_2  \left(\frac{k}{2n_1}  +\phi \right) \right]}   \cos   2\pi %n_3 \psi   =    \bar b_{k,l}(\phi)  \cos    2\pi n_3 \psi \\
%       \beta^{I}_{k,l} &=   \bar b_{k,l}(0)  \cos    2\pi n_3 \psi
 % \end{array}
  % \right.
%   \end{equation}
We can express   $D_1$  as a product of alternants in the  special case  where $2\pi n_3\psi = \pi/4$ and $  \phi = 0$. In this case
 equations \eqref{eq5} become 
 
 \begin{equation}\label{pisurquatre}
\left\{
\begin{array}{c}
 \alpha^I_{k,l}  =(-1)^{l+1}
   \frac{ \sin  \pi  \frac{n_3  k}{n_1}    }{\sqrt{2}\sin   \pi\frac{n_2 k}{n_1} }
   \left(  \cos\frac{\pi n_3l}{n_2}+  \sin\frac{\pi n_3 l}{n_2}\right)\\
\beta^I_{k,l} = (-1)^{l+1}  \frac{ \sin   \pi  \frac{n_3  k}{n_1}  }
   {\sqrt{2}\sin  \pi\frac{n_2 k}{n_1} } \left(  \cos\frac{\pi n_3l}{n_2}-  \sin\frac{\pi n_3 l}{n_2}\right) \\
  
\alpha^{II}_{k,l}  =(- 1)^{l+1} \frac{ \sin   \frac{\pi n_3 l}{n_2} }{\sqrt{2}\sin   \frac{\pi n_2 k}{n_1}}   
 \left( \cos\frac{\pi n_3 k}{n_1} +\sin   \frac{\pi n_3 k}{n_1}   \right) \\
  \beta^{II}_{k,l} = (- 1)^{l+1} \frac{ \cos     \frac{\pi n_3 l}{n_2} }{\sqrt{2}\sin  \frac{\pi n_2 k}{n_1}  }  \left( \cos\frac{\pi n_3 k}{n_1} +\sin   \frac{\pi n_3 k}{n_1}   \right)\end{array}
   \right.
   \end{equation}

% For ulterior use we will define  respectively \\
%  for points $(k;l)\in \Delta^I$ : 
%$$ a^{I}_{k,l} := (-1)^{l+1}  \frac{ \sin \left( \pi  \frac{n_3  k}{n_1}   \right)}
%   {\sin   \left(  \pi\frac{n_2 k}{n_1}\right) }\sin \left(  2\pi n_3 \left(  \frac{l}{n_2} +\psi\right) \right)$$
%\noindent and   for points $(k;l)\in \Delta^{II}$ : 
% $$ a^{II}_{k,l} = (- 1)^{l+1} \frac{\sin \left(  2\pi n_3  \left( \frac{k}{2n_1} +\phi +\psi\right) \right) }{\sin \left(  2\pi n_2  \left(\frac{k}{2n_1}  +\phi \right) \right)}\sin \left(  2\pi n_3  \left(   \frac{l}{n_2}\right) \right)$$

%Using Laplace expansion  for determinants,  we obtain, according to the parity of 
%$m$:
%$$D_1 = \sum_{\sigma \in\frak{S}_{2p}}\epsilon(\sigma)
%\left(\prod_{k=1}^{p}  \beta_{\sigma(p+k)  } \right)  \zeta ( a_{\sigma (1)}, \cdots, a_{\sigma (p)})
% \zeta ( a_{\sigma (p+1)}, \cdots, a_{\sigma (2p)})$$
 %or 
% $$D_1 = \sum_{\sigma \in\frak{S}_{2p+1}}\epsilon(\sigma)
%\left(\prod_{k=2}^{p}  \beta_{\sigma(p+k)  } \right)  \zeta ( a_{\sigma (1)}, \cdots, a_{\sigma (p+1)})
 %\zeta ( a_{\sigma (p+2)}, \cdots, a_{\sigma (2p+1)})$$
%\begin{equation}
%\epsilon(\sigma)  = \left( 
%\begin{array}{cccc}
%1, &2,&\cdots,&n \\
%\sigma(1), &\sigma(2),&\cdots,&\sigma(n)
%\end{array}
%\right)
%\end{equation}
%where $\zeta(a_1,\cdots, a_n) $ is the Sylvester notation for the difference product : $ \prod_{i<j} ( a_i -a_j )$
We compare the values of the $\alpha_P $  and $\alpha_Q $   ( resp.  $\beta_P$ and $\beta_Q $) for pairs of nodes as formed in section \ref{symetriesdesnoeuds}.
 More generally, let $\sigma$ be anyone of the following symmetries on the set of nodes $\Delta$ :
  \begin{center}
 \begin{tabular}{ccc}
 %Symmetries & on &  $\Delta^I$\\
 $\sigma_h^I :$ & $ (k,l) \mapsto$ &  $ (n_1-k,l)$\\
% $ \sigma_v^I:$ & $ (k,l) \mapsto$ &  $ (k,2n_2-l)$\\
%  $ \sigma_{v_1}^I:$& $ (k,l) \mapsto$ & $  (k, n_2 - l)$\\
%$ \sigma_{v_2}^I :$&$  (k,l) \mapsto$&$ (k, 3n_2 - l)$\\
     $ \tau^I : $&$  (k,l) \mapsto$&$ (k,l+n_2)$\\
 \end{tabular}
 \begin{tabular}{ccc}
 %Symmetries & on &  $\Delta^{II}$\\
 $\sigma_h^{II }:$ & $ (k,l) \mapsto$ &  $ (k,n_2-l)$\\
% $ \sigma_v^{II}:$ & $ (k,l) \mapsto$ &  $(2n_1-k, l)$\\
%  $ \sigma_{v_1}^{II}:$& $ (k,l) \mapsto$ & $(n_1-k, l)$\\
%$ \sigma_{v_2}^{II} :$&$  (k,l) \mapsto$&$ (3n_1-k,   l)$\\
     $ \tau^{II} : $&$  (k,l) \mapsto$&$ (k+n_1, l)$\\
 \end{tabular}
 \end{center}
  For  the sake of completeness we give the action of   symmetries  $\sigma$ on coefficients $\alpha$ or $\beta$ of 
  equations \eqref{pisurquatre} -defined by 
 $\sigma (\alpha_P) := \alpha_{\sigma(P)}$ where   $ P\in \Delta$. A straightforward  
  computation shows that 
 \begin{lem}\label{symm}  Let us assume  $n_2$ is odd   and $n_3$ is even.  Then 
 for values $\alpha$ and $\beta$ given in equations \eqref{pisurquatre} we have 
   \begin{center}
\begin{tabular}{c|c|c|}
 {\bf I} & $\alpha^I_{k,l}$ & $\beta^{I}_{k,l}$ \\
 \hline
 $\sigma_h^I$ & $  -\alpha^I_{k,l}$ &  $ -\beta^{I}_{k,l}$\\
 \hline 
% $ \sigma_v^I$ & $\beta^{I}_{k,l}$ &   $\alpha^I_{k,l}$\\
%   \hline
%%    $ \sigma_{v_1}^I$ & $\beta^{I}_{k,l}$ &   $\alpha^{I}_{k,l}$\\
%   \hline
%     $ \sigma_{v_2}^I$ &  $\beta^{I}_{k,l}$ &   $\alpha^{I}_{k,l}$\\
%   \hline
    $ \tau^I$ &  $-\alpha^{I}_{k,l}$ &   $-\beta^{I}_{k,l}$\\
   \hline

 \end{tabular}
 \qquad 
 \begin{tabular}{c|c|c|}
{\bf II} & $\alpha^{II}_{k,l}$ & $\beta^{II}_{k,l}$ \\
 \hline
 $\sigma_h^{II}$ & $ \alpha^{II}_{k,l}$ &  $ -\beta^{II}_{k,l}$\\
 \hline 
% $ \sigma_v^{II}$ & $ -\tilde\alpha^{II}_{k,l}$ &   $-\tilde\beta^{II}_{k,l}$\\
%   \hline
%   $ \sigma_{v_1}^{II}$ & $ \tilde\alpha^{II}_{k,l}$ &   $\tilde\beta^{II}_{k,l}$\\
%   \hline
%   $ \sigma_{v_2}^{II}$ & $ \tilde\alpha^{II}_{k,l}$ &   $\tilde\beta^{II}_{k,l}$\\
%   \hline
   $ \tau^{II}$ & $ -\alpha^{II}_{kl}$ &   $-\beta^{II}_{k,l}$\\
   \hline
 \end{tabular}
 \end{center}
 %where $\tilde\alpha^{II}_{kl}  =(- 1)^{l+1} \frac{ \sin   \frac{2\pi n_3 l}{n_2} }{\sqrt{2}\sin   \frac{\pi n_2 k}{n_1}}   
 %\left( \cos\frac{\pi n_3 k}{n_1} -\sin   \frac{\pi n_3 k}{n_1}   \right) ,$\\ and 
% $ \tilde\beta^{II}_{k,l} = (- 1)^{l+1} \frac{ \cos     \frac{2\pi n_3 l}{n_2} }{\sqrt{2}\sin  \frac{\pi n_2 k}{n_1}  }  \left( \cos\frac{\pi n_3 k}{n_1} -\sin   \frac{\pi n_3 k}{n_1}   \right) $
\end{lem}
 We deduce from Lemma \ref{symm}  that the values of the  $\alpha$ (resp.  $\beta$)  of two nodes that are coupled according 
 to Section \ref{coupling} are equal or opposite    (resp.  are opposite).
 With the help of this remark, we can find a simple expression for $D_0$ :\\
   let us regroup rows of the determinant $D_0$   by corresponding pairs of nodes as described in Subsection \ref{symetriesdesnoeuds} .  
We obtain two  such rows :
\begin{equation}\label{det1}
 \biggl|
\begin{array}{lllllll}
 1, &\alpha_{kl}^2,& \cdots,&\alpha_{kl}^{2p-2},& 
     \beta_{kl},& \beta_{kl}\alpha_{kl}^2,& \cdots\\
      1, &\alpha_{kl}^2,& \cdots,&\alpha_{kl}^{2p-2},& 
     -\beta_{kl},& -\beta_{kl}\alpha_{kl}^2,& \cdots\\
 \end{array}
      \biggr|    \end{equation}
  which- by elementary operations-  can be replaced by the following two rows in the determinant :  
  \begin{equation}\label{det2}
 \pm 4\beta_{kl} \biggl| 
      \begin{array}{cccccccc}
 1, &\alpha_{kl}^2,& \cdots,&\alpha_{kl}^{2p-2},& 
     0,& 0,& 0,&\cdots\\
       0,& \cdots ,& 0 ,& 0 & 
     1,&  \alpha_{kl}^2,& \cdots &\alpha_{kl}^{2p-2}\\
     \end{array}  
      \biggr|
  \end{equation}
  
  Doing this for all $p= n_1n_2 -(n_1+n_2)/2$  pairs-  parametrized by the two parallelograms $\Pi ^I \cup \Pi^{II}$- the determinant reduces to  the form 
    $\left|  
   \begin{array}{cc}
   A & 0 \\
   0 & A 
   \end{array}
   \right|$
   where $A$ is the Vandermonde matrix of size $p$-where  $m=2p$ :
     
      \begin{equation}\label{det3}
   A:= \biggl| 
   1,  \alpha_{kl}^2,  \cdots,\alpha_{kl}^{2p-2} \\
   \biggr|_{(k,l)\in \Pi^I \cup  \Pi^{II}}
   \end{equation}
\hfill \qed \\
\subsection{The lowest order term  $D_0$ of the Wronskian is not zero}\label{sub3}
We will show first that $D_0$ is not zero for phases  $\psi$ and $\phi$ of Lemma \ref{dezero}. 
 Since the Wronskian of the nodal curve $D$ is a polynomial in $r$ and non trivially nul -since $D_0 \not = 0$, this will show the non-nullity of $D$ for almost any rational $r$.
\begin{lem}\label{dezero}  Let  $n_1, n_2, n_3$  be relatively prime numbers such that 
$n_1$ and $n_2$ are odd primes  and  $n_3$   is even; suppose also that :

\begin{equation}\label{mod}
\left\{\begin{array}{c}
 n_2  \equiv 1\quad {\textrm mod\  } n_1 \\
 n_3  \equiv 2 \quad {\textrm mod\ } n_1 \\
   n_3  \equiv 2  \quad {\textrm mod\ } n_2
 \end{array}
 \right.
 \end{equation}
 Then the Wronskian   $D_0$ defined in \eqref{dzero} is not zero for values of $\alpha$'s and $\beta$'s defined 
 by equations \eqref{pisurquatre}.

     \end{lem}
      Let us notice   that the first congruence is true by Proposition \ref{residu} 
and  that the last two congruences are deduced from  the chinese remainder theorem.\\

    \noindent {\it Proof :} 
 %%%%%%%%%%%%%
Let us  check one by one that all the   factors of $D_0$ in  equation \eqref{dzero}  are not zero.
 \begin{itemize}
 \item 
Notice that  the coefficients $c_n$ defined in equations \eqref{loworder} are positive  ; indeed, by Euler's sine product formula 
 (cf.  for example  \cite{SZ}):
$$ \frac{u}{\sin u } = \frac{1}{\prod_{n=1}^\infty \left(  1- \frac{u^2}{\pi^2 n^2}  \right)},$$
Expanding the RHS, we obtain a series with even powers and positive coefficients. Hence all the even derivatives of $\frac{u}{\sin u }$ at zero are strictly positive. Postcomposing with  the $n$-th power 
 function  whose derivatives are all non-negatives at 1 and applying Faa di Bruno's expansion for  the derivatives of
  composed functions (cf.  for instance \cite{F}),  we show that   
  $$ \forall n\in \mathbb{N} \qquad c_n>0.$$
It is also clear that for almost any $\psi$ or $\phi$, the coefficients  $\alpha_{kl}$ are nonzero.
 \item  We check that  $\alpha^I_{kl} \not = 0$    since $n_1\ndiv k, n_2\ndiv l $ and $n_2$ and $n_3$ are relatively prime.
 We also check  that $\alpha^{II}_{kl} \not = 0$    since $n_1\ndiv k, n_2\ndiv l $ and $n_1$ and $n_2$ are odd   primes.
 \end{itemize}
 It remains to show that $D_1 \not = 0$:
 \begin{itemize}
 \item  We check that  $\beta^I_{kl} \not = 0$    since $n_1\ndiv k $ and $n_1$ and $n_3$ are coprime, 
 and that $\beta^{II}_{kl} \not = 0$    since $n_2 \ndiv l $ and $n_2$ and $n_3$ are coprime.
   \item   Let us prove that the difference-product  $\zeta(\alpha_j^2)$  defined in  Lemma \ref{dezerosept} is not  zero, i.e.  let us prove that the $\alpha^2_{j}$'s, where 
  $ j\in \Pi^I \cup  \Pi^{II}$,  are two by two distinct.\\
   Suppose on the contrary that there are distinct  $j=(k,l), j'=(k',l') \in \Pi^I \cup \Pi^{II}$  such that $\alpha_j = \pm \alpha_{j'}$. We will  prove that necessarily $j=j'$.\\
%  The coefficients  $\alpha_{j} $  can be simplified if computed in the cyclotomic field $\mathbb{Q}_{n_1n_28n_3}$.\\
%Consider   $\zeta := e^{\frac{ \pi i k}{n_1}}$,  $\xi:=   e^{\frac{ 2\pi i l}{n_2}}$, 
 % $ c:=   e^{\frac{\pi i}{8n_3}}$-where  $0<k< n_1, 0<l< n_2$-  which are respectively of order  resp. $n_1$ and $n_2$ ; hence 
 % $\zeta^{2n_1} = 1$ and    $\xi^{n_2} = 1$.\\
% We introduce the following  roots of unity :
% $$ \zeta =e^{i\frac{\pi k}{n_1}}, \xi = e^{i\frac{2\pi l}{n_2}}$$
%Since $n_1$ and $n_2$ are prime there  exist $s,r \in \mathbb{N}$ such that 
 %$$ \zeta^r =e^{i\frac{\pi k'}{n_1}}, \xi^s = e^{i\frac{2\pi l}{n_2}}$$
 
% Moreover, in each of the three cases that has to be considered, non trivial solutions will be sought 
%among  solutions $\zeta, \xi,\zeta^r,\xi^s$ whose  respective arguments 
%$\frac{\pi k}{n_1}, \frac{2\pi l}{n_2} ,\frac{\pi kr}{n_1},\frac{2\pi ls}{n_2} $ correspond to integer points 
%$(k,l), (kr,ls) $  that lie either in   $\Pi^I$ or $\Pi^{II}$ according to the studied cases.
Let us  examine successively the three possible cases :
  \begin{enumerate}
  \item
  Suppose $\alpha^I_{k,l} = \pm\alpha^I_{k',l'}$ i.e.
   \begin{equation}\label{cyclo1}
  \frac{ \sin  \pi  \frac{n_3  k}{n_1}    }{ \sin   \pi\frac{n_2 k}{n_1} }
   \left(  \cos\frac{\pi n_3l}{n_2}+  \sin\frac{\pi n_3 l}{n_2}\right) =
   \frac{ \sin  \pi  \frac{n_3  k'}{n_1}    }{ \sin   \pi\frac{n_2 k'}{n_1} }
   \left(  \cos\frac{\pi n_3l'}{n_2}+  \sin\frac{\pi n_3 l'}{n_2}\right)
   \end{equation}
   Using  congruences of \eqref{dezero}  and rearranging terms we obtain:
   \begin{equation}\label{cyclo2}
  \frac{ \sin  \pi  \frac{2 k}{n_1}    }{ \sin   \pi\frac{ k}{n_1} }\frac{ \sin   \pi\frac{ k'}{n_1} }{ \sin  \pi  \frac{2  k'}{n_1}    }
  =\frac{ \cos \pi  \frac{ k}{n_1}    }{ \cos  \pi\frac{ k'}{n_1} }
 = \pm \frac{ \cos\frac{2\pi l'}{n_2}+  \sin\frac{2\pi  l'}{n_2} }{     \cos\frac{2\pi l}{n_2}+  \sin\frac{2\pi l}{n_2}    }
   \end{equation} 
   Notice first that the RHS (resp. LHS)  is in the cyclotomic field $\mathbb{Q}_{2n_1}$ (in $\mathbb{Q}_{2n_2}$).
  Hence both  LHS and RHS are in $\mathbb{Q}_{2n_1}\cap \mathbb{Q}_{2n_2}=\mathbb{Q}$
   since $n_1$ and $n_2$ are prime.  
   Let us write the LHS in terms of   $\zeta := e^{\frac{ i\pi  }{n_1}}$ for some $ 0<k< \frac{n_1}{2} $; 
   Suppose $k\not = k'$. Without loss of generality we may suppose that  $k'> k$ and obtain  that for some  nonzero rational $q$ :
   \begin{equation}\label{ramadan}
   q\zeta^{2k'} - \zeta^{k'+k}- \zeta^{k'-k}+q=0.
   \end{equation} 
   The  cyclotomic polynomial $\frac{x^{n_1}+1}{x+1}$ divides   polynomial \eqref{ramadan}. Hence 
   $2k' > n_1 -1$ or $k'\geq \frac{n_1}{2}$. But then 
  $\zeta^{k'}$ doesn't belong to $\Pi^I$.\\
  Consequently, we must have  $k=k'$.\\
Thus the RHS equals $\pm 1$, from which  we deduce easily  that $l=l'$.
   % $\xi:=   e^{\frac{ 2\pi i l}{n_2}}$, 
     Hence   $\alpha_i \neq \pm\alpha_j$ whenever $i,j\in \Pi^{I}, i\neq j$.
    
   \item    
     Suppose $\alpha^{II}_{k,l} = \pm\alpha^{II}_{k',l'}$ i.e.
   \begin{equation}\label{cyclo1}
  \frac{ \sin   \frac{\pi n_3 l}{n_2} }{ \sin   \frac{\pi n_2 k}{n_1}}   
 \left( \cos\frac{\pi n_3 k}{n_1} +\sin   \frac{\pi n_3 k}{n_1}   \right) =
  \pm \frac{ \sin   \frac{\pi n_3 l'}{n_2} }{ \sin   \frac{\pi n_2 k'}{n_1}}   
 \left( \cos\frac{\pi n_3 k'}{n_1} +\sin   \frac{\pi n_3 k'}{n_1}   \right)
    \end{equation}
   Using  congruences of \eqref{dezero}  and rearranging terms we obtain:
      \begin{equation}\label{cyclo4}
  \frac{  \sin   \frac{\pi  k'}{n_1} \left( \cos\frac{2\pi  k}{n_1} +\sin   \frac{2\pi  k}{n_1}   \right)  }
  { \sin   \frac{\pi  k}{n_1}\left( \cos\frac{2\pi  k'}{n_1} +\sin   \frac{2\pi  k'}{n_1}   \right) } =
   \pm\frac{ \sin   \frac{2\pi  l'}{n_2} }{\sin   \frac{2\pi l}{n_2} }
   \end{equation} 
The first step is similar to the former case :  we  write the RHS in terms of  $\xi =e^{\frac{2\pi i l}{n_2}}$ 
and from  $\mathbb{Q}_{4n_1}\cap \mathbb{Q}_{n_2}=\mathbb{Q}$, we  deduce that $l=l'$.\\
 We derive  that 
  \begin{equation}\label{cyclo5}
  \frac{  \sin   \frac{\pi  k'}{n_1}}{\sin ( \frac{2\pi  k'}{n_1} +\frac{\pi}{4} )}  
 =\pm \frac {  \sin   \frac{\pi  k}{n_1}}{\sin ( \frac{2\pi  k}{n_1} +\frac{\pi}{4} )}.
   \end{equation} 
Plugging  $\xi = e^{i  \frac{\pi  k}{n_1}}, c = e^{i\frac{\pi}{4}},  \xi^s = e^{i  \frac{\pi  k'}{n_1}}$
in equation \eqref{cyclo5}, we obtain :
    \begin{equation}\label{cyclo7}
    (  \xi^2 - 1 ) (\xi^{4s} +i) -\sigma\xi^{s-1}(   \xi^{2s} - 1) (\xi^4 +i) =0, \quad \sigma =\pm 1.
 \end{equation} 

Hence $\xi $ is a root of a polynomial  \eqref{cyclo7}  with coefficients in 
$\mathbb{Q}(i)$.  But  the cyclotomic polynomial is irreductible in 
$\mathbb{Q}(i)$. Indeed the only quadratic extension of $\mathbb{Q}$ that lies in $\mathbb{Q}_{2n_1}$ is 
$\mathbb{Q}\left(\sqrt{(-1)^{\frac{n_1 -1}{2}}n_1}\right)$ hence 
$\mathbb{Q}(i)\cap \mathbb{Q}_{2n_1} = \mathbb{Q} $ ( cf. \cite{W} ). 
If the polynomial  \eqref{cyclo7}  is not trivially zero ( i.e $s=1, \sigma =1$)
then it  is a multiple 
of the cyclotomic polynomial, and  all the  primitive roots of unity are roots of polynomial
\eqref{cyclo7}. In particular $\xi^{-1} = \bar\xi$ is also a root. 
 We  combine  the conjugate of  equation \eqref{cyclo7} for $\xi$ and equation \eqref{cyclo7}  for $\xi^{-1}$;
letting  $\eta := \xi^4 $,  we derive that 
 \begin{equation}\label{cyclo8}
\frac{\eta^{s} - i}{\eta -i } =\frac{\eta^{s}+i}{\eta +i }.
 \end{equation} 
This implies immediatly that $\xi =1  $ or $\xi^s =1$.\\
 Finally   the polynomial \eqref{cyclo7} is trivially zero  i.e. $s=1$; thus $k=k'$.
Consequently  $\alpha_j \neq \pm\alpha_{j'}$ where $j,j'\in \Pi^{II}, j\neq j'$.
   \item    
    Suppose $\alpha^{I}_{k,l} = \pm\alpha^{II}_{k',l'}$ i.e.
   \begin{equation}\label{cyclo1}
 \frac{ \sin  \pi  \frac{n_3  k}{n_1}    }{ \sin   \pi\frac{n_2 k}{n_1} }
   \left(  \cos\frac{\pi n_3l}{n_2}+  \sin\frac{\pi n_3 l}{n_2}\right)=
  \pm \frac{ \sin   \frac{\pi n_3 l'}{n_2} }{ \sin   \frac{\pi n_2 k'}{n_1}}   
 \left( \cos\frac{\pi n_3 k'}{n_1} +\sin   \frac{\pi n_3 k'}{n_1}   \right)
    \end{equation}
   Using  congruences of \eqref{dezero}  and rearranging terms we obtain:
     \begin{equation}\label{cyclo1}
 \frac{ \sin  \pi  \frac{2  k}{n_1}    }{ \sin   \pi\frac{ k}{n_1} }
   \left(  \cos\frac{2\pi l}{n_2}+  \sin\frac{2\pi l}{n_2}\right)=
  \pm \frac{ \sin   \frac{2\pi  l'}{n_2} }{ \sin   \frac{\pi  k'}{n_1}}   
 \left( \cos\frac{2\pi  k'}{n_1} +\sin   \frac{2\pi  k'}{n_1}   \right)
    \end{equation}

      \begin{equation}\label{cyclo4}
 \frac{ \sin  \pi  \frac{2  k}{n_1}    }{ \sin   \pi\frac{ k}{n_1} }
 \frac{ \sin   \frac{\pi  k'}{n_1}}{ \left( \cos\frac{2\pi  k'}{n_1} +\sin   \frac{2\pi  k'}{n_1}   \right)}  =
   \pm \frac{\sin   \frac{2\pi  l'}{n_2}}{ \left(  \cos\frac{2\pi l}{n_2}+  \sin\frac{2\pi l}{n_2}\right)}
   \end{equation} 
Similarly the   RHS is rational;  we write it in terms of  $\xi =e^{\frac{2\pi i l'}{n_2}}$ and 
$\xi^s = e^{\frac{2\pi i l}{n_2}}$ for some $s$.
We obtain  
\begin{equation}\label{cyclo9}
\xi^{2s}  + q(1+i)\xi^{s+1}+q(1-i)\xi^{s-1}  -1=0.
\end{equation}
As before  the cyclotomic polynomial divides  polynomial \eqref{cyclo9} and any primitive root 
is a root.   In particular $\xi^{-1}$ is also a root ; summing the two equations we deduce that there is no solution with $\xi\neq 1$.
  \end{enumerate}
  Finally from the study of the  three cases, we deduce that no two  $\alpha_i$  for distinct $i$'s  are equal up to a sign for $i\in \Pi^I \cup \Pi^{II}$.
   \end{itemize}
  Whence  $D_1 \neq 0$, and Lemma \ref{dezero}  follows.   \hfill \qed\\

\subsection{$\mathbb{Q}$-linear dependence of nodal curves coordinates}
To finish the proof of Proposition \ref{proposition1}, it suffices to notice first that from subsection \ref{sub2}
the  Wronskian $D$ is a polynomial w.r.t. to the rational number $r =\frac{n_3}{n_2}$.
The coefficients of this polymonials are trigonometric expressions in $n_1, n_2, n_3$ but 
because of the congruences, these coefficients do not actually  depend on $n_3$.  Furthermore  the norm of   these  coefficients are uniformely bounded 
 from below and above by positive constants that depend only on $n_1$ and $n_2$. As we can choose  $n_3$, hence $r$ as large as we wish
we see that $D(r)$ is non-zero unless all coefficients are zero 
(we show that the coefficient of the monomial of highest order must be
zero and induction proves that all coefficients must be zero). Thus
    the coefficient of the monomial of lowest order in $r$ 
should be zero which  is not the case by Lemma \ref{dezero}. This shows Proposition \ref{proposition1}.\\
Finally we apply the infinitesimal deformation  criterion of Lemma \ref{courbe} and show that :
\begin{cor}\label{linindlis}
Any knot admits a family  of  Lissajous diagrams  such that the parameters of the nodes together with  one :\\
 $$\{ \eta^1(\epsilon), \cdots, \eta^{n_L}(\epsilon), 1 \}$$
are  $\mathbb{Q}$-linear independent   for infinitely many $\epsilon$.
\end{cor}

\section{Proof of the Theorem }
 
  We will prove following proposition, whence the theorem.
  \begin{prop}\label{termine}
  For any knot  $K$ there exist integers  $n_1,n_2,n_3,n_4$ 
  with $n_1,n_2$ odd primes and small positive numbers $\phi_0$ ,$\psi_0$
  such that for any $\phi, \psi$ ( and for almost any $\tau$) such that $|\phi|< \phi_0$, $|\psi -\frac{\pi}{4}|\leq \psi_0$, the curve defined by 
 \begin{equation}\label{plf}
 \gamma :
\left(  
\begin{array}{ll}
[0,1]  &\longrightarrow \mathbb{R}^3 \\
t  &\mapsto \left( \cos\left(2\pi n_1 t\right ) ,   \cos\left(2\pi n_2\left( t +\phi\right)\right ) +\epsilon 
\cos\left(2\pi n_3\left((t+\phi +\psi\right) \right),  \cos\left( 2\pi n_4\left( t+\tau \right) \right)   \right) 
\end{array}
\right)  
\end{equation} 
  is isotopic to $K$.
  \end{prop}
  We define a knot  $K$  by a knot shadow $D$ and  a {\it crossing sign function}
  $\alpha :\mathcal{N} \longrightarrow \{\pm 1\}$ which is  a  function from the set of nodes of 
  the shadow $D$.
   Each node $P_i \in \mathcal{N}$ ,$ i=1,\cdots, n_D$  of the shadow 
  is parametrized  by two parameters
   $(t_i,s_i), i=1,\cdots, n_D $.
 Our goal is to   construct a height function $z$ 
  such that       $sign\left( z\left(t_i\right) -z\left(s_i\right) \right) = \alpha(i) , \ i=1,\cdots, n $.\\

 From Proposition \ref{proposition1} there is a smooth knot 
 with  Lissajous shadow $  L(n_1,n_2,\phi)$ with  prime frequencies that is isotopic 
 to $K$.
  Then by  Corollary \ref{linindlis}  there is an $\epsilon$-deformation of $  L(n_1,n_2,\phi)$ :
   \begin{equation*}\label{plfbis}
L_\epsilon(n_1,n_2,n_3,\phi,\psi ):
\left(  
\begin{array}{ll}
[0,1]  &\longrightarrow \mathbb{R}^2 \\
t  &\mapsto \left( \cos\left(2\pi n_1 t\right ) ,   \cos\left(2\pi n_2\left( t +\phi\right)\right ) +\epsilon 
\cos\left(2\pi n_3\left(t+\phi +\psi\right) \right) \right) 
\end{array}
\right)  
\end{equation*} 
 such that the node parameters are  rationally linear independent. \\
 Let us prove now that we can find a frequency $n_4$ and  a  phase $\tau $ - and hence 
 %From Lemma \ref{hauteurlissajous}
   a height function    $\cos\left( 2\pi n_4\left( t+\tau \right) \right)$-   such that the  curve 
 defined in  Proposition \ref{termine}
is isotopic to $K$. This will prove that knot  $K$ is isotopic to a Fourier knot of type $(1,1,2)$.
 
\subsection{Kronecker's  theorem}
We will need the following direct  consequence of   Kronecker's theorem (cf. for instance \cite{HW}).
 \begin{lem}\label{kronembourg} Let $g : \mathbb{R} \longrightarrow J$ be a  1- periodic continuous function.\\
Let $t_1,\cdots, t_N, 1$ be $N$ real numbers  in $I:=[0,1] $ that are lin. ind. over $\mathbb{Q}$; 
and  let $\{v_i\}_{ i = 1, \cdots , N}$,   $N$ real values that lie in the image domain of $g$; then 
$$\forall \epsilon \ \forall  n_0\  \exists n>n_0 \quad |g(nt_i) - v_i | \leq \epsilon, \ i = 1, \cdots , N $$
Furthermore, we can choose $n$ relatively prime to   a given fixed  prime  integer $n_1$.
\end{lem}
\noindent{\it Proof}:
let $s_i\in I$ such that $g(s_i) = v_i \  i=1,\cdots, N $. 
We apply Kronecker's theorem (cf. \cite{HW}, Theorem 442) : for any $\epsilon$ and any $N$ 
there is an integer $n>n_0$  such that   
$$|n t_i -s_i  -p_i| \leq \epsilon\quad i =1,\cdots, N$$
where $  p_i\in \mathbb{Z}, i =1,\cdots, N  $.\\
Moreover  a modification  of    (cf. \cite{HW}, Theorem 200) allows us to show that 
 we can choose integer $n$ to be prime rel.  to $n_1$.\\
  Suppose  on the contrary that $n$ is a multiple of $n_1$; we define the vector $v := \left( t_1,\cdots, t_N \right)$
and consider the following sequence of vectors $v_j :=  l_j v\in [0,1]^N$ {\it  modulo one } where 
the integers  $l_j$ are chosen by induction so that  $l_{j+1} - l_{k}$ and  $n_1$  are 
 relatively prime $\forall k\leq j$.\\
 Divide the unit cube into subcubes with edges of length $1/Q$ for some $ Q>1$. For $ j=Q^N+1$ two of the $Q^N$ vectors $v_{i_1}, v_{i_2}$
  lie in the same subcube;  hence 
  $$|(l_{i_1}-l_{i_2}) t_i   -q_i| \leq \frac{\sqrt{2}}{Q^N}\leq \epsilon\qquad  \forall  i  = 1,\cdots, N$$
 where $q_i$ are integers . Let $m =(l_{i_1}-l_{i_2})$  and let $n' = m+n$ then $n'$ and $n_1$ are relatively prime and 
 $$|n' t_i -s_i  -p_i-q_i| \leq  2\epsilon \quad i =1,\cdots, N.$$
 \hfill \qed \\
The next lemma establishes the existence of a cosine function with the prescribed crossing signs.
\begin{lem}\label{hauteurlissajous} Let $D$ a be a Lissajous shadow  and let $\{ \left( s_i,t_i\right) \}$ 
be the set of $n_D$  nodes- double points-  distinct   parameters of $D$. 
Suppose that $\{ t_1, \cdots, t_{n_D}, 1\}$ are rationally lin. independent. 
Then for any  crossing  sign fonction $\alpha$ there exist a frequency $n_4$  and a phase $\tau$ such that the  function 
 \begin{equation}\label{definition de la quatrieme fonction trigonometrique}
 z\left(t\right) = \cos\left(2\pi n_4\left(t +\tau\right)\right)
 \end{equation}
 satisfies $sign\left( z\left(t_i\right) -z\left(s_i\right) \right) = \alpha(i) , \forall \ i=1,\cdots, n_D $.
\end{lem}
\noindent{\it  Proof} :  
 We first choose for each node one of the two parameters, say $t_1,\cdots, t_{n_D}$. It follows
 from Lemma \ref{kronembourg} that, for any $\epsilon$, there exists an integer $n$ with $gcd(n,n_1)=1$
 and 
 $|\cos(n t_i) -\alpha(i)|\leq \epsilon$.
  A node  of $D$ is parametrized by a $t_i$ and by the other parameter
 $s_i = \pm t_i + \frac{k}{n_1} \   i=1,\cdots, n_D,  0<k<n_1-1.$\\
 Suppose that $\alpha (i) =1$;  then there is a sequence of integers $n$ such that 
$\lim \cos(2\pi nt_i) = 1 \ \forall i =1,\cdots,n_D$. \\
We claim that for $n$ large enough 
$ \cos(2\pi t_i) > \cos(2\pi ns_i) \   \ \forall i =1,\cdots,n_D$.
Indeed 
 $$ \cos(2\pi ns_i) =   \cos(2\pi(\pm t_i))\cos\left(2\pi n   \frac{k}{n_1}\right) - \sin\left(2\pi n(\pm t_i\right)
 \sin\left(2\pi n \frac{k}{n_1}  \right) $$
The second term is close to zero when $n$ is large and the first  term is strictly less than  $|\cos\left(2\pi    \frac{k}{n_1}\right)|$
since $gcd(n,n_1)=1$. In conclusion whenever $\cos(2\pi nt_i)$ is close to $1$, $\cos(2\pi nt_i)-\cos(2\pi ns_i) >0$ and the strand corresponding 
to $z_n(t_i)$ is above the  strand corresponding to $z_n(s_i)$. The proof is identical if $\alpha(i)=-1$. We choose then $n$
such that $\cos(2\pi n(t_i)$ is close enough to $-1$.
\hfill \qed \\
In particular any knot is a (1,1,2) Fourier knot of the form 
 \begin{equation*}\label{plfbis}
\left(  
\begin{array}{ll}
[0,1]  &\longrightarrow \mathbb{R}^3 \\
t  &\mapsto \left( \cos2\pi n_1 t ,   \cos2\pi n_2\left( t +\phi\right) +\epsilon 
\cos\left(2\pi n_3t+\frac{\pi}{4}\right), \cos2\pi n_4 t\right)
\end{array}
\right)  
\end{equation*} 
for   positive numbers  $\phi$ and $\epsilon$ that are small enough.
  
\footnotesize{ 
Universit\'e F. Rabelais, LMPT  UMR 7750 CNRS   Tours, France,\\ marc.soret@lmpt.univ-tours.fr,
 Marina.Ville@lmpt.univ-tours.fr}
\end{document}